\documentclass{article}
\usepackage{titlesec}
\usepackage{titling}
\titleformat{\chapter}{\normalfont\huge}{\thechapter.}{20pt}{\huge\it}
\renewcommand{\paragraph}{\subsection}
\usepackage[utf8]{inputenc}
\usepackage[T1]{fontenc}
\usepackage{amsmath,amssymb,amsthm,mathrsfs,calligra,amsfonts}
\usepackage{enumerate}
\usepackage[framemethod=tikz]{mdframed}
\usepackage{array}
\usepackage{lipsum}
\usepackage{mathtools}
\usepackage{scrextend}
\usepackage{tikz-cd}
\usepackage{cite}
\usepackage{pdfrender,xcolor}
\usepackage{lmodern}
\usepackage[colorlinks]{hyperref}
\usepackage[noabbrev,capitalize]{cleveref}
\theoremstyle{plain}
\newtheorem{theorem}{Theorem}[section]   
\newtheorem{corollary}[theorem]{Corollary}
\newtheorem{lemma}[theorem]{Lemma}

\newtheorem{proposition}[theorem]{Proposition}
\hypersetup{
    colorlinks=true,       
    linkcolor=blue,        
    citecolor=green!70!black, 
    filecolor=magenta,     
    urlcolor=cyan          
}
\theoremstyle{definition}
\newtheorem{definition}[theorem]{Definition}

\newtheorem{example}[theorem]{Example}

\theoremstyle{remark}
\newtheorem{remark}[theorem]{Remark}
\titleformat*{\subsection}{\itshape\mdseries}
\crefformat{section}{\S#2#1#3}

\newcommand{\mr}{\mathrm}

\newcommand{\Z}{\mathbf Z}

\newcommand{\C}{\mathbf C}

\newcommand{\brak}{[\![\hbar]\!]}
\newcommand{\cbrak}{(\!(\hbar)\!)}

\newcommand{\2}{\vspace{2mm}}

\newcommand*{\shom}{\mathscr{H}\text{\kern -3pt {\it{om}}}\,}
\newcommand*{\ext}{\mathscr{E}\text{\kern -1.5pt {\it{xt}}}\,}
\makeatletter
\newcommand{\xdashrightarrow}[2][]{\ext@arrow 0359\rightarrowfill@@{#1}{#2}}
\newcommand{\xdashleftarrow}[2][]{\ext@arrow 3095\leftarrowfill@@{#1}{#2}}
\newcommand{\xdashleftrightarrow}[2][]{\ext@arrow 3359\leftrightarrowfill@@{#1}{#2}}
\def\rightarrowfill@@{\arrowfill@@\relax\relbar\rightarrow}
\def\leftarrowfill@@{\arrowfill@@\leftarrow\relbar\relax}
\def\leftrightarrowfill@@{\arrowfill@@\leftarrow\relbar\rightarrow}
\def\arrowfill@@#1#2#3#4{%
  $\m@th\thickmuskip0mu\medmuskip\thickmuskip\thinmuskip\thickmuskip
   \relax#4#1
   \xleaders\hbox{$#4#2$}\hfill
   #3$%
}
\makeatother
\usepackage{geometry}
\usepackage[parfill]{parskip}    
\usepackage{graphicx}
\usepackage{epstopdf}
\title{\vspace{-1cm}Deformations of cohomology classes and spectral sequences}
\author{Borislav Mladenov}
\AtEndDocument{\bigskip{\footnotesize%
  \text{Institute of Mathematics, Academia Sinica, Taipei, Taiwan} \\
  \href{mailto:mladenov@as.edu.tw}{\textit{mladenov@as.edu.tw}}
  }}
\predate{}
\date{}
\postdate{}
\begin{document}
\maketitle
\nocite{*}

\begin{abstract}
We study spectral sequences as an obstruction theory to deforming cohomology classes via the Rees construction. A filtered complex is interpreted as an equivariant family over the affine line, in which higher spectral sequence differentials arise as obstructions to lifting cohomology classes through infinitesimal neighborhoods of the special fibre. This leads to various geometric and module-theoretic degeneration criteria. We give various applications to $\mr{Ext}$ spectral sequences for complex Lagrangians and the formality results of Melani-Rubi\'o amongst others.\\
$\hspace*{5mm}$For double complexes, we revisit the Khovanov-Qi decomposition on the Rees side and introduce a two-parameter Rees construction over the affine plane. The higher differentials of the two associated spectral sequences correspond to torsion modules supported on the thickened coordinate axes. We obtain degeneration criteria in terms of these modules and an intrinsic characterisation of the $\partial\bar\partial$-lemma.
\end{abstract}
\hypersetup{linkcolor=black} 
\tableofcontents

\hypersetup{linkcolor=blue} 
\section{Introduction}

\subsection{Spectral sequences via the Rees functor}
$\hspace*{3mm}$A filtered complex contains two complexes: its underlying complex and the associated graded complex representing the classical limit of the filtration. The spectral sequence attached to the filtration measures the failure of cohomology to commute with passage to the associated graded. The purpose of this note is to make this failure geometric.\\
$\hspace*{3mm}$The basic device is the Rees construction. Let
$$
(\mr{C},\mr{F})
$$
be a decreasingly filtered complex over a field $k$, its Rees complex is
$$
\mr{R}_\mr{F}(\mr{C})
=
\bigoplus_p \mr{F}^p\mr{C}\,\hbar^{-p}.
$$
It is naturally a graded complex over $k[\hbar]$, or equivalently, a
$\mathbf G_m$-equivariant complex on the affine line. Its special fibre at
$\hbar=0$ is
$$
\mr{R}_\mr{F}(\mr{C})\otimes_{k[\hbar]}k
\simeq
\operatorname{Gr}_\mr{F}(\mr{C}),
$$
whereas at $\hbar=1$ one recovers the original
complex $\mr{C}$. Thus, we may regard a filtered complex as a one-parameter
degeneration from $\mr{C}$ to $\operatorname{Gr}_\mr{F}(\mr{C})$ and the associated spectral sequence measures the failure of cohomology to commute with the classical limit $\hbar =0$: $$\mr{R}_\mr{F}\left(\mr{H}(\mr{C})\right)/\hbar \stackrel{?}{\simeq}\mr{H}\left(\mr{R}_\mr{F}(\mr{C})/\hbar\right).$$ 
$\hspace*{3mm}$In order to work entirely on the Rees side, we express the pages of the spectral sequence $(\mr{E}_r,\mr{d}_r)$ associated to $(\mr{C,F})$ in terms of the Rees object. The Rees complex of a filtered complex is naturally bigraded: $$\mr{R}_\mr{F}^n(\mr{C}) = \bigoplus_{p\in \Z}\mr{R}_\mr{F}^{p,n-p}(\mr{C}), \text{\;\; where\;\; } \mr{R}_\mr{F}^{p,n-p}(\mr{C})=\mr{F}^p\mr{C}^n\hbar^{-p}.$$ Compatibility of the differential with the filtration implies cohomology is bigraded and we define $$\mr{H}^{p,q}(\mr{R}_\mr{F}(\mr{C}))=\mr{H}^{p+q}\left(\mr{R}_\mr{F}^{p,\bullet-p}(\mr{C})\right).$$ Multiplication by $\hbar$ has bidegree $(-1,1)$. In particular, the cokernel and its cohomology are naturally bigraded too.
\2
 \begin{proposition}[\cref{pages}]
 Let $(\mr{C,F})$ be a filtered complex with an associated spectral sequence $(\mr{E}_r,\mr{d}_r)$. For all $r\ge 1$ and all $p,q \in \Z$, there are natural isomorphisms
 $$\mr{E}_r^{p,q} \cong \mr{im}\left(\cdot \hbar^{r-1}: \mr{H}^{p,q}\left(\mr{R_\mr{F}(C)}/\hbar^r\right)\to \mr{H}^{p-r+1,q+r-1}\left(\mr{R_\mr{F}(C)}/\hbar^r\right)\right).$$
 \end{proposition}
 $\hspace*{3mm}$Similarly, corresponding to these isomorphisms, there is a Rees description of the differentials which goes as follows:
\2 
 \begin{proposition}[\cref{diffs}]
 Let $(\mr{C,F})$ be a filtered complex with an associated spectral sequence $(\mr{E}_r,\mr{d}_r)$. For all $r\ge 1$ and all $p,q \in \Z$ there are commutative diagrams:
 \[
\begin{tikzcd}
\mr{E}_r^{p,q} \arrow[r, "\mr{d}_{r}"] \arrow["\cong",d]
& \mr{E}_r^{p+r,q-r+1}\arrow[d,"\cong"]   \\
\mr{im}\left(\cdot \hbar^{r-1}|_{\mr{H}^{p,q}\left(\mr{R_\mr{F}(C)}/\hbar^r\right)}\right)\arrow[r, "\delta_{r,r}"]
& \mr{im}\left(\cdot \hbar^{r-1}|_{\mr{H}^{p+r,q-r+1}\left(\mr{R_\mr{F}(C)}/\hbar^r\right)}\right).
\end{tikzcd}
\]
where $\delta_{r,r}$ is the connecting morphism associated to the short exact sequence: $$0\xrightarrow{} \mr{R}_\mr{F}(\mr{C})/\hbar^r \xrightarrow{\cdot \hbar^r}  \mr{R}_\mr{F}(\mr{C})/\hbar^{2r} \xrightarrow{} \mr{R}_\mr{F}(\mr{C})/\hbar^r \xrightarrow{} 0.$$
\end{proposition}
$\hspace*{3mm}$These elementary observations are the organising principle of the paper. We interpret spectral sequence differentials as obstruction classes for extending cohomology classes away from the special fibre, identify degeneration with flatness properties of Rees cohomology, and compare a number of apparently different deformation-theoretic phenomena. The resulting picture includes filtered derived categories, Hard Lefschetz, deformation quantisation, decompositions of double complexes and the $\partial\bar\partial$-lemma, deformation of cohomology classes in the sense of Wavrik, and the $T^1$-lifting of Ran and Kawamata.

\subsection{Flatness and degeneration}

$\hspace*{3mm}$The first topic we address is that degeneration is a flatness statement on the Rees side. Indeed, this is the analogue of Deligne's criterion: $\mr{E}_1$-degeneration is equivalent to the strictness of the differential. 
\2
\begin{theorem}[\cref{e1flatness}]
Let $(\mr{C,F},\mr{d})$ be a biregular filtered complex. Its spectral sequence degenerates at $\mr{E}_1$ if and only if $\mr{H}(\mr{R_\mr{F}(C)})$ is a flat $k[\hbar]$-module.
\end{theorem}
The exact sequences $$0 \longrightarrow \mr{H}^{p,q}(\mr{R_\mr{F}(C)})/\hbar \longrightarrow \mr{E}_1^{p,q} \longrightarrow \mr{Tor}^{k[\hbar]}_1(\mr{H}^{p,q+1}(\mr{R_\mr{F}(C)}),k)\longrightarrow 0.$$
show that failure of an $\mr{E}_1$-class to lift to Rees cohomology is measured by $\hbar$-torsion in the adjacent degree. Hence, as the (total) cohomology is graded over $k[\hbar]$, having no $\hbar$-torsion is equivalent to being flat, and so we obtain:
\2
\begin{corollary}[\cref{lifts}]
Let $(\mr{C,F},\mr{d})$ be a biregular filtered complex. Its spectral sequence degenerates at $\mr{E}_1$ if and only if every $\mr{E}_1$-class extends to  $\mr{H}\left(\mr{R_\mr{F}(C)}\right)$.
\end{corollary}

\subsection{Obstructions to infinitesimal lifting problems}

$\hspace*{3mm}$In the discussion above, we are extending classes on the central fibre to global sections, skipping all intermediate finite order lifts. The second theme is addressing these infinitesimal lifting problems by realising the differentials
\[
\mr{d}_r:\mr{E}_r^{p,q}\longrightarrow \mr{E}_r^{p+r,q-r+1}.
\]
as obstructions to the finite order lifting problems. In addition, this will show that the biregularity condition on the filtration is unnecessarily strong and can be replaced with a much weaker one.
$\hspace*{3mm}$Let
\[
\mr{R}_{\mr{F},N}(\mr{C}):=\mr{R}_\mr{F}(\mr{C})/\hbar^{N+1}.
\]
be the truncated Rees complex associated to the filtered complex $\mr{(C,F)}$.
Given a class on the central fibre
\[
x_0\in
\mr{H}^{p+q}(\operatorname{Gr}_\mr{F}\mr{C})
\] one can ask if it extends to a class in $\mr{R}_N$, i.e. to a class over the $N^\text{th}$-infinitesimal neighbourhood. Under the identification of the $E_1$-page with
central fibre cohomology, our statement reads as follows:
\2
\begin{theorem}[\cref{infinitesimal}]
Let $(\mr{C,F},\mr{d})$ be a filtered complex with associated spectral sequence $\left(\mr{E}_r,\mr{d}_r\right)$. A class $x_0$ on the central fibre of the Rees complex, that is $x_0\in \mr{E}_1$, extends to the $N^{\text{th}}$-infinitesimal neighbourhood if and only if $$\mr{d}_1\left([x_0]\right)=0,\;\; \mr{d}_2\left([x_0]\right)=0,\cdots, \mr{d}_N\left([x_0]\right)=0.$$
Thus the higher differentials are the successive obstruction classes for extending a
cohomology class to the infinitesimal neighbourhoods of the special fibre.
\end{theorem}
$\hspace*{3mm}$In particular, we see that flatness implies $\mr{E}_1$-degeneration under no additional assumptions on the complex or the filtration:
\2
\begin{theorem}[\cref{flate1}]
Let $(\mr{C,F},\mr{d})$ be a filtered complex. Assume that $\mr{H}(\mr{R}_\mr{F}(\mr{C}))$ is flat over $k[\hbar]$. Then the spectral sequence degenerates at $\mr{E}_1$.
\end{theorem}
$\hspace*{3mm}$For the opposite direction, examining the proofs below, one sees that the only obstruction is infinite $\hbar$-divisible torsion. Hence we introduce a separatedness condition on the cohomology ruling out such divisibility:
\2
\begin{theorem}[\cref{separatede1}]
Let $(\mr{C,F},\mr{d})$ be a filtered complex and suppose its spectral sequence collapses at $\mr{E}_1$. Assume that the Rees cohomology $\mr{H}(\mr{R}_\mr{F}(\mr{C)})$ is $\hbar$-adically separated: $$\bigcap_{N\ge 1} \hbar^{N}\mr{H}(\mr{R}_\mr{F}(\mr{C)})=0.$$ Then $\mr{H}(\mr{R}_\mr{F}(\mr{C)})$ is flat.
\end{theorem}
$\hspace*{3mm}$Finally, we show one can work over the formal neighbourhood of the origin instead.
\2
\begin{theorem}[\cref{formale1}]
Let $(\mr{C,F},\mr{d})$ be a filtered complex. Assume that $\mr{H}\left(\widehat{\mr{R}}_\mr{F}(\mr{C})\right)$ is flat over $k\brak$; equivalently, every $\mr{E}_1$-class has a formal lift. Then the spectral sequence degenerates at $\mr{E}_1$.\\
Conversely, if the spectral sequence degenerates at $\mr{E}_1$ and $\mr{H}\left(\widehat{\mr{R}}_\mr{F}(\mr{C})\right)$ satisfies the separateness condition: $$\bigcap_{N\ge 1} \hbar^{N}\mr{H}\left(\widehat{\mr{R}}_\mr{F}(\mr{C)}\right)=0,$$ then $\mr{H}\left(\widehat{\mr{R}}_\mr{F}(\mr{C})\right)$ is flat and every $\mr{E}_1$-class has a formal lift.
\end{theorem}
\2
\begin{remark}
Note that lifting (independently) to every finite thickening, constructing a compatible inverse system of lifts, producing a class in the (completed) Rees complex, and comparing cohomology of the completed Rees complex with the completion of Rees cohomology are distinct problems and the first is controlled purely by the spectral sequence, while the latter require completeness, separatedness, or perfectness hypotheses and do not depend on the spectral sequence.
\end{remark}
\2
\begin{remark}
The same obstruction mechanism appears in classical deformation theory.
Wavrik's obstruction theory for deforming a class
\[
\xi_0\in \mr{H}^q(\mr{X}_0,\mathscr F_0)
\]
in a flat proper analytic family is obtained from the successive extensions
\[
A/\mathfrak m^{n+1}\to A/\mathfrak m^n
\]
and takes values in
\[
\mr{H}^{q+1}(\mr{X}	_0,\mathscr F_0)\otimes
\mathfrak m^n/\mathfrak m^{n+1},
\]
see \cite{WavrikDeformingCohomologyClasses}. So, over one parameter, the Rees construction identifies Wavrik's successive obstructions with the differentials of the spectral sequence.
\end{remark}
\subsection{Multiplication and degeneration}

$\hspace*{3mm}$ Up next, we use that the Rees functor is monoidal and hence, for a filtered dg algebra $\mr{(A,F,d)}$, it gives a filtered dg algebra $\mr{R}_\mr{F}(\mr{A})$. In \cref{mult}, we observe that using the multiplicative structure, it is enough to lift a set of generator for the algebra $\mr{H}(\mr{G}_\mr{F}(\mr{A}))$ to $\mr{H}(\mr{R}_\mr{F}(\mr{A}))$ in order to prove the latter is flat over $k[\hbar]$. In particular, we apply our obstruction theory and the resulting spectral sequence degeneration criteria to the $\mr{Ext}$ spectral sequence of a compact K\"ahler Lagrangian $\mr{L}$ in a holomorphic symplectic $\mr{X}$. 
\2
\begin{theorem}[\cref{e2lag}]
Let $\mr{X}$ be a holomorphic symplectic variety and let $i:\mr{L}\xhookrightarrow{} \mr{X}$ be a compact K\"ahler Lagrangian whose Hodge algebra $$\bigoplus_{p,q}\mr{H}^p(\Omega_\mr{L}^q)$$ is generated by classes of holomorphic degree $q\le 1$. Suppose $\mathscr{L}$ is a line bundle on $\mr{L}$. Then the $\mr{Ext}$ spectral sequence $$\mr{H}^p(\Omega_\mr{L}^q) \Rightarrow \mr{Ext}^{p+q}(i_*\mathscr{L},i_*\mathscr{L})$$ collapses at $\mr{E}_2$ if and only if $\mr{d}_2^{p,1}=0$ for all $p\ge 0$.
\end{theorem}
$\hspace*{3mm}$A particularly simple case is that of abelian Lagrangians:
\2
\begin{corollary}[\cref{abeliancase}]
Let $\mr{X}$ be a holomorphic symplectic variety and let $i:\mr{L}\xhookrightarrow{} \mr{X}$ be an abelian Lagrangian. Let $\mathscr{L}$ be a line bundle on $\mr{L}$. Then the $\mr{Ext}$ spectral sequence $$\mr{H}^p(\Omega_\mr{L}^q) \Rightarrow \mr{Ext}^{p+q}(i_*\mathscr{L},i_*\mathscr{L})$$ collapses at $\mr{E}_2$ if and only if $\mr{d}_2^{0,1}=0$. This hold for example when $\mathscr{L}=\mathscr{O}_\mr{L}$.
\end{corollary}
\2
\begin{remark}
These results are similar to the main theorems of \cite{MladenovDegeneration}. The main difference is that we impose conditions on the cohomology of the Lagrangian and require more than just $\mr{d}_2(\omega)=0$, but we no longer need $\mathscr{L}$ to extend to first order in every direction.
\end{remark}
\subsection{Flatness via Hard Lefschetz}

$\hspace*{3mm}$It becomes a natural question how to prove Rees flatness without computing the spectral sequence
differentials. When our complexes arise from geometry, one of the main mechanisms is Hard Lefschetz.\\
$\hspace*{3mm}$Suppose a perfect complex $\mr{K}$ over a principal ideal domain $R$ carries a Serre duality pairing
\[
\mr{K}\simeq \mr{R}\operatorname{Hom}(\mr{K},k[\hbar])[-2n]
\]
and a Lefschetz operator
\[
\mr{L}:\mr{H}^k(\mr{K})\longrightarrow \mr{H}^{k+2}(\mr{K})
\]
satisfying
\[
\mr{L}^{n-k}:\mr{H}^k(\mr{K})\xrightarrow{\sim}\mr{H}^{2n-k}(\mr{K}).
\]
Duality shifts torsion by one complementary degree, while Hard Lefschetz
reflects it across the middle dimension. Combining the two, we obtain that
any nonzero torsion class gives rise to nonzero torsion classes in all degrees. Since our complex is bounded, we see immediately that no torsion can occur.
\2
\begin{theorem}[\cref{HL-duality-flatness}]
Let $R$ be a principal ideal domain and let
$\mr{K}\in \mathrm{Perf}(R)$. Assume that
\[
\mr{K}\simeq \mr{R}\operatorname{Hom}_R(\mr{K},R)[-2n]
\]
and that a degree-two operator $\mr{L}$ satisfies Hard Lefschetz
\[
\mr{L}^{\,n-i}:\mr{H}^i(\mr{K})\xrightarrow{\sim}\mr{H}^{2n-i}(\mr{K})
\]
for all $i\leq n$. Then every $\mr{H}^i(\mr{K})$ is torsion-free. Hence every
$\mr{H}^i(\mr{K})$ is a finite free $R$-module.
\end{theorem}

$\hspace*{3mm}$This turns the construction of a Lefschetz operator into a degeneration
criterion. The same philosophy underlies the classical K\"ahler examples. For a
variation of Hodge structure, Simpson identifies the filtered de Rham complex with a
$\lambda$-connection family whose central fiber is the corresponding Higgs
complex. Strictness of the Hodge filtration corresponds to flatness of the Rees
family, while the K\"ahler identities provide the Lefschetz structure which
forces that flatness.

\subsection{Formality of homotopy algebras}

$\hspace*{3mm}$Another application concerns formality of homotopy algebras. Recall that a homotopy algebra over a dg operad is \emph{formal} if it is connected to its cohomology, endowed with the induced strict algebra structure, by a zigzag of quasi-isomorphisms. In a minimal homotopy model, higher operations measure the possible failure of formality. The filtration by arity on its deformation complex leads to an operadic cohomology spectral sequence. Building on Manetti's criterion for differential graded Lie algebras, Melani and Rubi\'o show that formality is characterised by degeneration of this spectral
sequence at $\mr E_2$, see \cite{ManettiFormality,MelaniRubio}.\\
$\hspace*{3mm}$The Rees construction explains this criterion in terms of lifting a distinguished Euler class. Let $\mathcal{P}$ be a dg Koszul operad and let $(\mr{A,Q})$ be a minimal homotopy algebra over $\mathcal{P}$. Let $\mathfrak{g}$ be the deformation complex of $\mr{A}$ filtered by weight. After d\'ecalage of the weight filtration, we can trivialise the complete Rees complex $\mathfrak{g}$ as $$\widehat{\mr{R}}_{\mr{Dec(F)}}(\mathfrak{g}) \simeq \mathfrak{g}\brak$$
and the Rees differential becomes
\[ \mr{D}_{\mr{Q},\hbar}=[\mr{Q}_\hbar,-],
\]
where $\mr{Q}_\hbar$ is the deformation to the normal cone given by 
\[
\mr{Q}_\hbar=\mr{q}_2+\hbar \mr{q}_3+\hbar^2 \mr{q}_4\cdots.
\]
A calculation shows
\[
\mr{D}_{\mr{Q},\hbar}(\mr{e_A})=\hbar\,\partial_\hbar\mr{Q}_\hbar.
\]
The cocycle $\partial_\hbar\mr{Q}_\hbar$ represents the Kaledin class. Hence, we have:
\2
\begin{proposition}[\cref{ek}]
Let $\mr{(A,Q)}$ be a minimal homotopy $\mathcal{P}$-algebra. The Kaledin class $\mr{K}_\mr{A}$ is the obstruction class to lifting the Euler class $\mr{e}_\mr{A}$ to a formal Rees cocycle.
\end{proposition}
\2
\begin{remark}
The successive finite order lifting obstructions are thus the truncated Kaledin classes. See also \cite{KaledinFormalityFamilies,EmprinKaledin}.
\end{remark}
$\hspace*{3mm}$Working on the Rees side, we obtain that the Kaledin class measures the failure of the Rees cohomology to be flat over $k\brak$:
\2
\begin{theorem}
Let $\mathcal{P}$ be a reduced dg Koszul operad of finite type. Let $\mr{(A,Q)}$ be a minimal homotopy $\mathcal{P}$-algebra such that $\mr{dim}_k(\mr{A)}<\infty$. The following are equivalent:
\begin{enumerate}
\item The Kaledin class $\mr{K}_\mr{A}$ vanishes.
\item The Rees cohomology $\mr{H}\left(\mathfrak{g}\brak,\mr{D}_{\mr{Q},\hbar}\right)$ is flat over $k\brak$.
\end{enumerate}
\end{theorem}
\2
\begin{remark}
Paired with the main theorem of \cite{MelaniRubio}, we deduce that $\mr{(A,Q)}$ is formal if and only if the Rees cohomology $\mr{H}\left(\mathfrak{g}\brak,\mr{D}_{\mr{Q},\hbar}\right)$ is flat.
\end{remark}

\subsection{Double complexes and the two-parameter Rees construction}

$\hspace*{3mm}$Next, we specialise to double complexes. Let
\[
(\mr{C}^{p,q},\partial,\bar\partial).
\]
be a double complex, its total complex carries two natural filtrations and we introduce a bi-Rees construction
\[
\mr{R}_{s,t}(\mr{C})
=
\bigoplus_{p,q}
(\mr{F}^p\operatorname{Tot}\mr{C}\cap\mr{G}^q\operatorname{Tot}\mr{C})\,s^{-p}t^{-q}.
\]
Using the bigrading, we see it is quasi-isomorphic to
\[
\left(
\operatorname{Tot}(\mr{C})\otimes k[s,t],
s\partial+t\bar\partial
\right).
\]
Specialising $t=1$, after forgetting the second grading, gives the ordinary
Rees complex for the column filtration; specialising $s=1$ gives the row
Rees complex. The specialisations at $s=0$ and $t=0$ are the
associated-gradeds, and at their intersection
\[
(s,t)=(0,0)
\]
we get the bigraded central fibre.\\
$\hspace*{3mm}$We start by interpreting the decomposition of bounded double complexes into indecomposables of
\cite{StelzigDoubleComplexes,KhovanovQiTricomplexes} on the (bi)-Rees side. Since the Rees functor is additive, it's enough to evaluate on the indecomposables.\\
$\hspace*{3mm}$For a one-parameter Rees complex, this is done is \cref{khoqi}. Simples and odd length zigzags correspond to rank 1 free Rees modules, squares give contractible complexes, and for even length zigzags, we get a torsion module or a contractible complex for horizontal and vertical orientations, respectively. The Khovanov-Qi decomposition is not canonical but is finer than the decomposition of the cohomology Rees module over the PID $k[\hbar]$. The latter, however, is canonical and exists for filtered complexes too. If one cares only about the spectral sequence, the cohomology Rees module is the natural object to consider. In particular, we immediately recover that degeneration for a spectral sequence associated to a double complex is equivalent to having no horizontal even zigzags.\\
$\hspace*{3mm}$As explained, the one-parameter Rees cohomology does not distinguish between simples and odd length zigzags. It also doesn't see the squares. On the other hand, recall a theorem of Stelzig that a double complex satisfies the $\partial\bar\partial$-lemma if and only if its decomposition contains only simples and squares. In order to give a geometric interpretation of this result which does not depend on choosing a decomposition of the double complex, we need to refine our approach and use the two-parameter Rees functor.\\
$\hspace*{3mm}$ Let $R=k[s,t]$ and $\mathfrak{m}=(s,t)$. The elementary pieces become:
\[
\begin{array}{c|c}
\text{indecomposable}
&
\text{two-parameter Rees module}
\\ \hline
\text{simple}
&
k[s,t]
\\
\text{square}
&
\mr{K}_R(s,t)
\\
\text{even zigzag}
&
\text{module on a thickened coordinate axis}
\\
\text{odd zigzag}
&
\mathfrak{m}^r \text{ or its derived dual},
\end{array}
\]
where $\mr{K}_R(s,t)$ is the Koszul complex.
Thus spectral sequence differentials appear in codimension one, while the failure of the two filtrations to be opposed is supported in
codimension two. This leads to a geometric characterisation of the $\partial\bar\partial$-lemma and the main result goes as follows:
\2
\begin{theorem}[\cref{ddctwo}]
Let $\mr{C}$ be a bounded double complex of finite-dimensional vector spaces and
let
\[
\mr{R}_{s,t}(\mr{C})
\]
be its two-parameter Rees complex over $R=k[s,t]$.
The following are equivalent.

\begin{enumerate}
\item $\mr{C}$ satisfies the $\partial\bar\partial$-lemma.

\item In the decomposition of $\mr{C}$, only simples and squares
occur.

\item The following intrinsic conditions hold:
\begin{enumerate}
\item[(a)] $\mr{H}\left(\mr{R}_{s,t}(\mr{C})\right)$ is locally free on
\[
\mathbf A^2\setminus\{0\};
\]
\item[(b)] for every $m$,
\[
\mr{H}^m\left(\mr{R}_{s,t}(\mr{C})\right)/H^0_{\mathfrak m}(\mr{H}^m\left(\mr{R}_{s,t}(\mr{C})\right))
\]
is reflexive;
\item[(c)]
\[
\mathfrak m\,H^0_{\mathfrak m}\left(\mr{H}^m\left(\mr{R}_{s,t}(\mr{C})\right)\right)=0
\]
for every $m$;
\item[(d)] Conditions (a)--(c) above hold and the remaining Postnikov classes vanish.
\end{enumerate}

\item Conditions (a)--(c) above hold and the hyper-Tor spectral
sequence
\[
\operatorname{Tor}^R_p\left(\mr{H}^q\left(\mr{R}_{s,t}(\mr{C})\right),k\right)
\Longrightarrow
\mr{H}^{q-p}\left(\left(\mr{R}_{s,t}(\mr{C})\right)\otimes_R^Lk\right)
\]
degenerates at $E_2$.
\end{enumerate}
\end{theorem}
We can express these conditions in terms of sheaves on the affine plane, making the conditions much more familiar to algebraic geometers. Let
\[
R=k[s,t],\qquad
\mr{M}^m=\mr{H}^m\left(\mr{R}_{s,t}(\operatorname{C})\right),\qquad
\mr{T}^m=H^0_{\mathfrak{m}}\left(\mr{M}^m\right).
\]
The
\begin{enumerate}
\item $\widetilde{\mr{M}^m}$ is locally free on
\[
\mathbf A^2\setminus\{0\};
\]
equivalently, both spectral sequences degenerate at $\mr{E}_1$.

\item the quotient
\[
\widetilde{\mr{M}^m}/\widetilde{\mr{T}^m}
\]
satisfies the Hartogs' property at the origin.

\item the torsion subsheaf $\widetilde{\mr{T}^m}$ is scheme-theoretically supported on
the reduced origin.
\item The singularities at the origin are split in the sense that the extension classes of free and skyscraper cohomology sheaves vanish.
\end{enumerate}
\2
\begin{remark} 
The first condition removes even zigzags; the second removes vertical odd zigzags;
the third then gets rid of the horizontal odd zigzags for $r\ge 2$. This leaves us with simples, squares and three-vertex odd zigzags. The latter is distinguished from a sum of a simple and a square precisely because it corresponds to a non-trivial $\mr{Ext}^2$-class between free and skyscraper cohomology modules, so we require these to vanish.\\
We also note that, for a complex satisfying the first three conditions, the last one is equivalent to formality of the complex.
\end{remark}
\subsection{Acknowledgements}
$\hspace*{3mm}$This note has been long overdue. Many thanks go to Richard Thomas who suggested the main philosophy developed here a few years ago, but the text was never in a state to be made public.\\
$\hspace*{3mm}$ The author is supported by a postdoctoral fellowship at the Institute of Mathematics, Academia Sinica.
\section{Filtrations and the Rees functor}\label{sec1}
$\hspace*{3mm}$We begin by recalling some basic notions on filtered and graded algebras and modules to fix notation and terminology.
\subsection{Filtrations}
\begin{definition}
\begin{enumerate}
\item A filtered vector space is a pair $(\mr{V},\mr{F})$, where $\mr{V}$ is a vector space and $\mr{F}$ is a $\Z$-indexed decreasing filtration on $\mr{V}$.
\item A morphism of filtered vector spaces $f:(\mr{V},\mr{F})\to (\mr{W},\mr{F})$ is a linear map $f:\mr{V}\to\mr{W}$ such that for any $p\in \Z$ we have $f(\mr{F}^p\mr{V})\subseteq \mr{F}^p\mr{W}$.
\item We write $\mr{Fil}(\mr{Vect})$ for the category of filtered vector spaces and their morphisms.
\end{enumerate}
\end{definition}
$\hspace*{3mm}$Given two filtered vector spaces $(\mr{V,F})$ and $(\mr{W,F})$, their tensor product $\mr{V}\otimes \mr{W}$ is naturally a filtered vector space via $$\mr{F}^k(\mr{V}\otimes\mr{W})\coloneq \bigoplus_{p+q=k}\mr{F}^p\mr{V}\otimes\mr{F}^q\mr{W}.$$ This gives $\mr{Fil}(\mr{Vect})$ a symmetric monoidal structure.
\2 
\begin{definition}
A filtered algebra is an algebra in $\mr{Fil}(\mr{Vect})$. Explicitly, it is a filtered vector space $(\mr{A},\mr{F})$ such that $\mr{A}$ is an algebra and the multiplication is a morphism of filtered vector spaces $(\mr{A\otimes A},\mr{F})\to(\mr{A},\mr{F})$.
\end{definition}
\2
\begin{definition}
Let $(\mr{A},\mr{F})$ be a filtered algebra. 
\begin{enumerate}
\item A filtered $\mr{A}$-module is a filtered vector space $(\mr{M},\mr{F})$ such that $\mr{M}$ is an $\mr{A}$-module and the multiplication $\mr{A}\otimes\mr{M}\to \mr{M}$ is a morphism of filtered vector spaces.
\item A morphism of filtered $\mr{A}$-modules $f:(\mr{M,F})\to (\mr{N,F})$ is a morphism of $\mr{A}$-modules $f:\mr{M}\to\mr{N}$ compatible with the filtrations, i.e. for any $p\in \Z$ we have $f(\mr{F}^p\mr{M})\subseteq \mr{F}^p\mr{N}$.
\item We denote the category of filtered $\mr{A}$-modules by $\mr{Fil}(\mr{Mod}(\mr{A}))$.
\end{enumerate}
\end{definition}

\subsection{Gradings}
\begin{definition} Let $\mr{V}$ be a vector space.
\begin{enumerate}
 \item A $\Z$-grading on the vector space $\mr{V}$ is a decomposition into subspaces $$\mr{V} = \bigoplus_{n \in \Z} \mr{V}_n.$$ A $\Z$-graded vector space is a vector space equipped with a $\Z$-grading.
\item A morphism of graded vector spaces $f:\mr{V}\to \mr{W}$ is a linear map $f:\mr{V}\to\mr{W}$ preserving the gradings, i.e. $f(\mr{V}_n)\subseteq\mr{W}_n$.
\item We write $\mr{Gr}(\mr{Vect})$ for the category of graded vector spaces.
\end{enumerate}
\end{definition}
$\hspace*{3mm}$Given graded vector spaces $$\mr{V} = \bigoplus_{n \in \Z} \mr{V}_n \text{ and } \mr{W} = \bigoplus_{n \in \Z} \mr{W}_n,$$ we get a grading on the tensor product $$\mr{V\otimes W}=\bigoplus_{n\in \Z}\bigoplus_{p+q=n}\mr{V}_p\otimes\mr{W}_q$$ which gives $\mr{Gr}(\mr{Vect})$ a symmetric monoidal structure. 
\2
\begin{definition}
A graded algebra is an algebra in $\mr{Gr}(\mr{Vect})$. Equivalently, it is a graded vector space $$\mr{A}=\bigoplus_{n\in \Z}\mr{A}_n$$ such that $\mr{A}$ is an algebra and the multiplication $$\mr{A}\otimes\mr{A} \to \mr{A}$$ is a morphism of graded vector spaces.
\end{definition}
\2
\begin{definition}
Let $\mr{A}$ be a graded algebra.
\begin{enumerate}
\item A graded $\mr{A}$-module is a graded vector space $$\mr{M}=\bigoplus_{n\in \Z}\mr{M}_n$$ such that $\mr{M}$ is an $\mr{A}$-module and the multiplication $$\mr{A}\otimes\mr{M}\to \mr{M}$$ is a morphism of graded vector spaces.
\item A morphism of graded $\mr{A}$-modules $f:\mr{M}\to \mr{N}$ is a morphism of $\mr{A}$-modules $f:\mr{M}\to\mr{N}$ compatible with the gradings, i.e. for any $n\in \Z$ we have $f(\mr{M}_n)\subseteq \mr{N}_n$. 
\item We write $\mr{Gr}(\mr{Mod}(\mr{A}))$ for the category of graded $\mr{A}$-modules and their morphisms.
\end{enumerate} 
\end{definition}
\2
\begin{remark}[Associated graded]
 For any $p\in \Z$, there is a functor $\mr{Gr}^p : \mr{Fil}(\mr{Vect})\to \mr{Vect}$ sending $$\mr{(V,F)} \mapsto \mr{F}^p\mr{V}/\mr{F}^{p+1}\mr{V}.$$ Piecing these together, we get an additive functor $\mr{Gr}:\mr{Fil}(\mr{Vect})\to \mr{Gr}(\mr{Vect})$ given by $$(\mr{V,F})\mapsto \bigoplus_{p\in \Z}\mr{F}^p\mr{V}/\mr{F}^{p+1}\mr{V}.$$ It is a symmetric monoidal functor and the natural morphism $$\mr{Gr}(\mr{V,F})\otimes \mr{Gr}(\mr{W,F})\to \mr{Gr}(\mr{V\otimes W,F})$$ is an isomorphism. For a filtered algebra $\mr{(A,F)}$, we get an induced functor $$\mr{Gr}:\mr{Fil}(\mr{Mod}(\mr{A}))\to \mr{Gr}(\mr{Mod}(\mr{Gr}\mr{A})).$$
\end{remark}
\2
\begin{remark}
We usually write $\mr{Gr_F(M)}$ for the graded module $\mr{Gr}(\mr{M,F})$.
\end{remark}

\subsection{The Rees construction}
$\hspace*{3mm}$ We are now ready to introduce the Rees functor. 
\2
\begin{definition}
Let $(\mr{A,F})$ be a filtered algebra. The Rees algebra $\mr{R}_\mr{F}(\mr{A})$ is the graded subring $$\mr{R}_\mr{F}(\mr{A}) = \bigoplus_{p\in \Z} \mr{F}^{p}\mr{A}\cdot\hbar^{-p} \subset \mr{A}\otimes k[\hbar,\hbar^{-1}].$$
Let $(\mr{M,F})$ be a filtered $(\mr{A,F})$-module. The Rees module of $(\mr{M,F})$ is the graded $\mr{R}_\mr{F}(\mr{A})$-module $$\mr{R}_\mr{F}(\mr{M})=\bigoplus_{p\in \Z} \mr{F}^{p}\mr{M}\cdot\hbar^{-p}\subset \mr{M}\otimes k[\hbar,\hbar^{-1}].$$ 
\end{definition}
\begin{remark}
If $f:(\mr{M,F})\to(\mr{N,F})$ is a filtered morphism of $\mr{(A,F)}$-modules, we get a morphism of graded $\mr{R}_\mr{F}(\mr{A})$-modules $\mr{R}(f):\mr{R}_\mr{F}(\mr{M})\to\mr{R}_\mr{F}(\mr{N})$, turning the Rees construction into a functor $$\mr{Fil}(\mr{Mod}(\mr{A,F}))\to\mr{Gr}(\mr{Mod}(\mr{R}_\mr{F}(\mr{A}))).$$ 
In particular, we have a Rees functor $$\mr{Fil}(\mr{Vect})\to \mr{Gr}(\mr{Mod}(k[\hbar])).$$
\end{remark}
\begin{proposition}
The Rees functor is fully faithful and its essential image consists of graded $\mr{R}_\mr{F}(\mr{A})$-modules with no $k[\hbar]$-torsion.
\end{proposition}
\begin{proof}
Let $g : \mr{R}(\mr{M,F})\to\mr{R}(\mr{N,F})$ be a morphism of graded $\mr{R}(\mr{A,F})$-modules. Hence $$g_p : \mr{F}^p\mr{M}\cdot\hbar^{-p} \to \mr{F}^p\mr{N}\cdot\hbar^{-p}$$ and since $g$ is compatible with the module structure, it follows that $g$ commutes with $\hbar$. Since $\varinjlim \mr{F}^p\mr{M}=\mr{M}$ and $\varinjlim \mr{F}^p\mr{N}=\mr{N}$, setting $$f \coloneqq \varinjlim g_p :\mr{M}\to\mr{N},$$ we get $\mr{R}_\mr{F}(f)=g$. This shows the functor is full. Similarly, if $f_1,f_2:\mr{M} \to \mr{N}$ are filtered morphisms such that $\mr{R}_\mr{F}(f_1)=\mr{R}_\mr{F}(f_2)$, then $f_1$ and $f_2$ agree on $\mr{F}^p\mr{M}$ for all $p \in \Z$ and so $f_1=f_2$.\\
$\hspace*{3mm}$To show that the essential image is precisely $k[\hbar]$-torsion free graded $\mr{R}_\mr{F}(\mr{A})$-modules, we let $$\widetilde{\mr{M}}=\bigoplus_{p\in \Z} \widetilde{\mr{M}}^p$$ be such module. By assumption, it has no $\hbar$-torsion, so $\hbar : \widetilde{\mr{M}} \to \widetilde{\mr{M}}(-1)$ is injective. This in turn is equivalent to $\widetilde{\mr{M}}^p\subseteq \widetilde{\mr{M}}^{p-1}$ for all $p\in \Z$. Now, define $$\mr{M} \coloneqq \varinjlim \widetilde{\mr{M}}^{p}.$$ This is a filtered $\mr{R}_\mr{F}(\mr{A})$-module and by construction $\mr{R_\mr{F}(M)}=\widetilde{\mr{M}}$.
\end{proof}
\2
\begin{definition}
Let $(\mr{A},\mr{F})$ be a filtered algebra.
\begin{enumerate}
\item A flat Rees module is a graded $\mr{R}_\mr{F}(\mr{A})$-module which is $k[\hbar]$-flat.
\item A morphism between flat Rees modules is flat if its kernel and cokernel are flat Rees modules.
\item A Rees complex is a complex of graded $\mr{R}_\mr{F}(\mr{A})$-modules. We call it flat if its cohomology modules are flat Rees modules.
\end{enumerate}  
\end{definition}
\2
\begin{remark}
By the previous proposition, the flat Rees modules are precisely these graded $\mr{R}_\mr{F}(\mr{A})$-modules which are in the image of the Rees functor.
\end{remark}
$\hspace*{3mm}$ The next result will be important for applications to the study of spectral sequences of filtered complexes.
\2
\begin{proposition}\label[proposition]{flatness}
Let $(\mr{A},\mr{F})$ be a filtered algebra. A complex of flat Rees modules is a flat Rees complex if and only if its differentials are flat (as morphisms between flat Rees modules).
\end{proposition}
\begin{proof}
Suppose $$\cdots \xrightarrow{\mr{d}^{i-1}} \mr{C}^i\xrightarrow{\mr{d}^i}\mr{C}^{i+1}\xrightarrow{}\cdots$$ is a bounded above complex of flat Rees modules. Let $$\mr{B}^i=\operatorname{im}\mr{d}^{i-1},\qquad \mr{Z}^i=\ker \mr{d}^i.$$ These are by hypothesis $k[\hbar]$-flat. We have a short exact sequence $$0\longrightarrow \mr{H}^i(\mr{C})
\longrightarrow
\operatorname{coker}(\mr{d}^{i-1})
\longrightarrow
\mr{B}^{i+1}
\longrightarrow0.
$$
If the Rees complex is flat, each cohomology is $k[\hbar]$-flat, hence so is the extension $\operatorname{coker}(\mr{d}^{i-1})$, so the differentials are flat.\\
$\hspace*{3mm}$ Conversely, if the latter have $k[\hbar]$-flat cokernels for all $i$, then the $k[\hbar]$-submodule $\mr{H}^i(\mr{C})$ is also $k[\hbar]$-flat. Hence the Rees complex is flat.
\end{proof}
$\hspace*{3mm}$The next result is a practical criterion for checking flatness of a morphism.
\2
\begin{proposition}\label[proposition]{flatmorphism}
A morphism $\varphi : \mr{M} \to \mr{N}$ of flat graded $k[\hbar]$-modules is flat if and only if it is $\hbar$-pure, i.e. $\varphi(\mr{M})\cap \hbar\mr{N}=\hbar\varphi(\mr{M})$.
\end{proposition}
\begin{proof}
A graded $k[\hbar]$-module is flat if and only if it has no $\hbar$-torsion. The $\hbar$-purity of $\varphi$ is equivalent to saying, for $[v]\in \operatorname{coker}(\varphi)$, if $\hbar\cdot [v]=0 \in \operatorname{coker}(\varphi)$, then $v\in \operatorname{im}(\varphi)$, i.e. $[v]=0 \in \operatorname{coker}(\varphi)$. In other words $\operatorname{coker}(\varphi)$ has no $\hbar$-torsion. Since it is graded, the claim follows.
\end{proof}
\2
\begin{proposition}\label[proposition]{basecase}
 We have the following isomorphisms:
 \begin{enumerate}
  \item $\mr{R}_\mr{F}(\mr{M})/(\hbar-1) \cong \mr{M}$ as (filtered) $(\mr{A,F})$-modules.
  \item $\mr{R}_\mr{F}(\mr{M})[\hbar^{-1}] \cong \mr{M} \otimes k[\hbar,\hbar^{-1}]$ as $\mr{A}\otimes k[\hbar,\hbar^{-1}]$-modules.
  \item $\mr{R}_\mr{F}(\mr{M})/\hbar \cong \mr{Gr_F}(\mr{M})$ as graded $\mr{Gr_F}(\mr{A})$-modules.
 \end{enumerate}
\end{proposition}
\begin{proof}
\begin{enumerate}
\item Consider the evaluation at $\hbar=1$ $$\mr{ev}_1 :\mr{R}_\mr{F}(\mr{M}) \to \mr{M}.$$It is surjective as the filtration is exhaustive. Let $$v = \sum_{p=a}^{b}v_p\cdot\hbar^{-p} \in \operatorname{ker}(\mr{ev}_1).$$ Then $$\mr{ev}_1(v) = \sum_{p=a}^{b}v_p$$ and define $$\mr{S}_r = \sum_{p=a}^{r}v_p= -\sum_{p=r+1}^{b}v_p \in \mr{F}^{p+1}\mr{M}.$$Hence, letting $$w= \sum_{r=a}^{b-1}S_r\hbar^{-(r+1)},$$ we have $v=(z-1)w$ which shows that $$(\hbar-1)\mr{R_\mr{F}(M)} = \operatorname{ker}(\mr{ev}_1).$$ We recover the filtration from the grading of the Rees module: $$\mr{F}^p\mr{M} = \operatorname{im}\left(\mr{R_\mr{F}(M)}_{-p} \to \mr{R_\mr{F}(M)}/(\hbar-1)\right).$$
\item This is by definition.
\item Indeed, the degree $-p$ piece of $\mr{R_\mr{F}(M)}$ is $$\mr{F}^p\mr{M}\cdot\hbar^{-p}.$$
Similarly, the $-p$ piece of $\hbar\mr{R_\mr{F}(M)}$ is $$\hbar\left(\mr{F}^{p+1}\mr{F}\cdot\hbar^{-(p+1)}\right)=\mr{F}^{p+1}\mr{M}\cdot\hbar^{-p}.$$ Hence $$\left(\mr{R_\mr{F}(M)}/\hbar\right)_{-p}=\frac{\mr{F}^p\mr{M}}{\mr{F}^{p+1}\mr{M}}\cdot \hbar^{-p}$$ and so we have $$\mr{R_\mr{F}(M)}/\hbar=\bigoplus_{p\in \Z}\frac{\mr{F}^p\mr{M}}{\mr{F}^{p+1}\mr{M}}\cdot\hbar^{-p}.$$ It is clear that this is compatible with the module structures: if $a \in \mr{F}^p\mr{A},\;\; v \in \mr{F}^q\mr{M}$, then $$av \in \mr{F}^{p+q}\mr{M}.$$ For the Rees modules $$\left(a\cdot\hbar^{-p}\right)\left(v\cdot\hbar^{-q}\right) = av\cdot \hbar^{-p-q}.$$ Modulo $\hbar$ this becomes precisely $[a]\cdot[v]=[av]$ in the associated gradeds.\qedhere
\end{enumerate}
\end{proof}
$\hspace*{3mm}$For some of the applications, we need to work with unbounded filtrations, so it will be useful to introduce a few variations of the Rees functor that are suited to the non-biregular setting.
\2
\begin{definition}[Pro and completed Rees functors]
Let $(\mr{A,F})$ be a filtered algebra. For $N\ge 0$, its $N$-truncated Rees algebra is $$\mr{R}_{\mr{F},N}\mr{(A)} = \mr{R}_\mr{F}\mr{(A)}/\hbar^{N+1}.$$
Let $\mr{(M,F)}$ be a filtered $\mr{(A,F)}$-module. Similarly, the $N$-truncated Rees module is $$\mr{R}_{\mr{F},N}\mr{(M)} = \mr{R}_\mr{F}\mr{(M)}/\hbar^{N+1}.$$ 
The pro-Rees algebra and module are defined as the pro-systems given by these truncations: $$\mathbf{R}_\mr{F}(\mr{A}) = \{\mr{R}_{\mr{F},N}\mr{(A)}\}_{N\ge 0} \text{ and } \mathbf{R}_\mr{F}(\mr{M}) = \{\mr{R}_{\mr{F},N}\mr{(M)}\}_{N\ge 0}.$$
The completed Rees algebra and module are the inverse limits $$\widehat{\mr{R}}_\mr{F}(\mr{A}) = \varprojlim_N \mr{R}_{\mr{F},N}\mr{(A)} \text{ and } \widehat{\mr{R}}_\mr{F}(\mr{M}) = \varprojlim_N \mr{R}_{\mr{F},N}\mr{(M)}.$$
\end{definition}
\section{Spectral sequences associated to filtrations}

$\hspace*{3mm}$We begin this section by briefly introducing spectral sequences associated to filtered differential modules following \cite{cartan1956homological}.
\subsection{Filtered differential modules}
$\hspace*{3mm}$Let $(\mr{A,F})$ be a filtered algebra and let $(\mr{M,F,d})$ be a filtered differential module over it. If we ask for a grading on $\mr{M}$ and $\mr{d}$ to have degree $1$ with respect to the grading, we get a filtered complex of $\mr{A}$-modules. Let $\mr{H(M)}$ be the cohomology of $(\mr{M,d})$. We have a natural filtration on it: $$\mr{F}^p\mr{H(M)}  = \mr{im}(\mr{H}(\mr{F}^p\mr{M}) \to \mr{H(M)}).$$ The components of the associated graded are given by:
\2
\begin{proposition}
 There are natural isomorphisms $$\mr{Gr}_{\mr{F}}^p(\mr{H(M)}) \simeq \mr{im}(\mr{H}(\mr{F}^p\mr{M}) \to \mr{H}(\mr{M}/\mr{F}^{p+1}\mr{M})).$$
\end{proposition}
$\hspace*{3mm}$For all $n\ge 1$, we define the space of $n$-cocycles $$\mr{Z}^p_n = \mr{im}(\mr{H}(\mr{F}^p\mr{M}/\mr{F}^{p+n}\mr{M}) \to \mr{H}(\mr{F}^p\mr{M}/\mr{F}^{p+1}\mr{M}))$$ and the space of $n$-coboundaries $$\mr{B}_n^p = \mr{im}(\mr{H}(\mr{F}^{p-n+1}\mr{M}/\mr{F}^{p}\mr{M}) \to \mr{H}(\mr{F}^p\mr{M}/\mr{F}^{p+1}\mr{M})).$$ Since we have $\mr{B}_n^p \subseteq \mr{Z}_n^p$, we can define $\mr{E}_n^p = \mr{Z}_n^p/\mr{B}_n^p.$ Using the natural isomorphisms $$\mr{Z}_n^p/\mr{Z}_{n+1}^p \simeq \mr{im}(\mr{H}(\mr{F}^p\mr{M}/\mr{F}^{p+n}\mr{M}) \to \mr{H}(\mr{F}^{p+1}\mr{M}/\mr{F}^{p+n+1}\mr{M}))  \simeq \mr{B}_{n+1}^{p+n}/\mr{B}_{n}^{p+n},$$ we get a map
\2
 \begin{equation*}\mr{d}_n^p :\mr{E}_n^p \xhookrightarrow{} \mr{Z}_n^p/\mr{Z}_{n+1}^p \simeq \mr{B}_{n+1}^{p+n}/\mr{B}_{n}^{p+n} \to \mr{E}_n^{p+n}.\end{equation*}
 
It follows that $$\mr{ker}\left(\mr{d}_n^p\right) = \mr{Z}_{n+1}^{p}/\mr{B}_{n}^{p} \text{\quad and \quad} \mr{im}\left(\mr{d}_n^p\right) = \mr{B}_{n+1}^{p+n}/\mr{B}_{n}^{p+n},$$ hence $$\mr{ker}\left(\mr{d}_n^p\right)/\mr{im}\left(\mr{d}_n^{p-n}\right) \simeq \mr{E}_{n+1}^{p}.$$ 
$\hspace*{3mm}$Define the space of $\infty$-cocycles by $$\mr{Z}^p_{\infty} = \mr{im}(\mr{H}(\mr{F}^p\mr{M}) \to \mr{H}(\mr{F}^p\mr{M}/\mr{F}^{p+1}\mr{M}))$$ and the space of $\infty$-coboundaries
$$\mr{B}^p_{\infty} = \mr{im}(\mr{H}(\mr{M}/\mr{F}^{p}\mr{M}) \to \mr{H}(\mr{F}^p\mr{M}/\mr{F}^{p+1}\mr{M})).$$ As above, we set $$\mr{E}^p_{\infty}  = \mr{Z}^p_{\infty}/\mr{B}^p_{\infty}.$$ There is a chain of inclusions $$\cdots \subseteq \mr{B}^p_{n} \subseteq \mr{B}^p_{n+1}\subseteq \cdots \mr{B}^p_{\infty} \subseteq \mr{Z}^p_{\infty} \subseteq \cdots \subseteq \mr{Z}^p_{n+1}\subseteq \mr{Z}^p_{n}\subseteq \cdots$$ and if we define $\mr{E}^p_{0} = \mr{Gr}^p_{\mr{F}}\mr{M}$ with $\mr{d}_0^p$ induced by $\mr{d} : \mr{M} \to \mr{M}$, we have:
\2
\begin{proposition}[\cite{cartan1956homological}]	
 For all $n\ge 0$, let $\mr{E}_n = \bigoplus_p \mr{E}_n^p$ and $\mr{d}_n = \bigoplus_p \mr{d}_n^p$. Then $\mr{d}_n$ is a differential of degree $n$ and $\mr{H}(\mr{E}_n,\mr{d}_n) \simeq \mr{E}_{n+1}$.
\end{proposition}
$\hspace*{3mm}$We would like to give a different description of the differentials that has an easy interpretation on the equivariant side.\\
The next proposition follows from the above definitions and the diagram:
\[
 \begin{tikzcd}[row sep = large]
  & \mr{H}(\mr{F}^{p}\mr{M}/\mr{F}^{p+n}\mr{M})\arrow[dr]\arrow[d]\\
  \mr{H}(\mr{F}^{p-n+1}\mr{M}/\mr{F}^{p}\mr{M}) \arrow[r]\arrow[ur] & \mr{H}(\mr{F}^{p}\mr{M}/\mr{F}^{p+1}\mr{M})\arrow[r] & \mr{H}(\mr{F}^{p-n+1}\mr{M}/\mr{F}^{p+1}\mr{M}).
 \end{tikzcd}
\]

\begin{proposition}\label[proposition]{ssiso}
 For any $1\le n \le \infty$, there are isomorphisms $$\mr{E}_n^p \simeq \mr{im}(\mr{H}(\mr{F}^p\mr{M}/\mr{F}^{p+n}\mr{M}) \to \mr{H}(\mr{F}^{p-n+1}\mr{M}/\mr{F}^{p+1}\mr{M})).$$
\end{proposition}
$\hspace*{3mm}$Consider now the diagram
\[
 \begin{tikzcd}
  0\arrow[r] & \mr{F}^{p+n}\mr{M}/\mr{F}^{p+2n}\mr{M} \arrow[r]\arrow[d] & \mr{F}^{p}\mr{M}/\mr{F}^{p+2n}\mr{M} \arrow[r]\arrow[d] & \mr{F}^{p}\mr{M}/\mr{F}^{p+n}\mr{M} \arrow[r]\arrow[d] & 0\\
  0\arrow[r] & \mr{F}^{p+1}\mr{M}/\mr{F}^{p+n+1}\mr{M} \arrow[r] & \mr{F}^{p-n+1}\mr{M}/\mr{F}^{p+n+1}\mr{M} \arrow[r] & \mr{F}^{p-n+1}\mr{M}/\mr{F}^{p+1}\mr{M} \arrow[r]& 0.
 \end{tikzcd}
\]
Taking cohomology, we see that the connecting homomorphism of the bottom sequence induces a homomorphism between the images of the vertical maps:
\[
 \begin{tikzcd}
\mr{H}(\mr{F}^{p}\mr{M}/\mr{F}^{p+n}\mr{M}) \arrow[r] \arrow[d]
& \mr{H}( \mr{F}^{p+n}\mr{M}/\mr{F}^{p+2n}\mr{M}) \arrow[d] \\
\mr{H}( \mr{F}^{p-n+1}\mr{M}/\mr{F}^{p+1}\mr{M}) \arrow[r]
& [] \mr{H}( \mr{F}^{p+1}\mr{M}/\mr{F}^{p+n+1}\mr{M}).
\end{tikzcd}
\]
Using $\cref{ssiso}$, the bottom horizontal map induces a map $\mr{E}_n^p \to \mr{E}_n^{p+n}$ which is precisely the differential.

\subsection{Deligne's \texorpdfstring{$\mr{E}_1$}--degeneration}

$\hspace*{3mm}$We have the following simple criterion which is surprisingly geometric as we will see later.
\2
\begin{proposition}[Deligne]\label[proposition]{e1}
The spectral sequence $(\mr{E}_r,\mr{d}_r)$ associated to a biregular filtered differential module $\mr{(M,F,d)}$ collapses at $\mr{E}_1$ if and only if the differential $\mr{d}$ is strict with respect to the filtration $\mr{F}$, in other words, $$\mr{d}(\mr{F}^p\mr{M})=\mr{d}(\mr{M})\cap\mr{F}^p\mr{M}.$$
\end{proposition}
\begin{proof}
Suppose $\mr{d}$ is strict. Take a cocycle $\bar{x} \in \mr{Gr_F}^p\mr{M}$, so $\mr{d}(x) \in \mr{F}^{p+1}\mr{M}$. By strictness of $\mr{d}$, we can find some $y\in \mr{F}^{p+1}\mr{M}$ with $$\mr{d}(x)=\mr{d}(y).$$ Then $x-y$ represents the same class as $x$ in $\mr{Gr_F}^p\mr{M}$, so $$\operatorname{ker}(\mr{Gr_F}\mr{d})= \mr{Gr_F}\left(\operatorname{ker}(\mr{d})\right).$$ Similarly, we have $$\operatorname{im}(\mr{Gr_F}\mr{d})= \mr{Gr_F}\left(\operatorname{im}(\mr{d})\right)$$ which implies that $$\mr{H}\left(\mr{Gr_F}\mr{M}\right) = \mr{Gr_F}\left(\operatorname{ker}(\mr{d})\right)/\mr{Gr_F}\left(\operatorname{im}(\mr{d})\right).$$ The last quotient is also isomorphic to $\mr{Gr_F}(\mr{HM})$, hence $$\mr{H}\left(\mr{Gr_F}\mr{M}\right) \simeq \mr{Gr_F}(\mr{HM})$$ and thus the spectral sequence collapses at $\mr{E}_1$.\\
Suppose the spectral sequence degenerates on $\mr{E}_1$, so $\mr{d}_r=0$ for all $r$. Define $x^{(1)}=x$ and suppose $$x^{(r)} \equiv x \pmod {\mr{F}^{p+1}\mr{M}}$$ is such that $\mr{d}\left(x^{(r)}\right) \in \mr{F}^{p+r}\mr{M}$. It defines a class in $\mr{E}_r$ and the vanishing $\mr{d}_r\left(x^{(r)}\right)=0$ implies that there is $y^{(r)} \in \mr{F}^{p+1}\mr{M}$ with $$\mr{d}\left(x^{(r)}-y^{(r)}\right) \in \mr{F}^{p+r+1}\mr{M}.$$ Now, we define $$x^{(r+1)} =x^{(r)}-y^{(r)}.$$ By finiteness of the filtration, we can find $r$ large enough such that $$\mr{d}\left(x^{(r)}\right)=0.$$ Letting $y= \mr{d}(x)=\mr{d}\left(x-x^{(r)}\right)$, we see $$y \in \mr{d}\left(\mr{F}^{p+1}\mr{M}\right),$$ so that $$\mr{d}(\mr{F}^{p+1}\mr{M})=\mr{d}(\mr{F}^p\mr{M})\cap\mr{F}^{p+1}\mr{M}.$$ Iterating the latter, we get $$\mr{d}(\mr{F}^{p+r}\mr{M})=\mr{d}(\mr{F}^p\mr{M})\cap\mr{F}^{p+r}\mr{M},$$ for all $r\ge 1$. This implies $\mr{d}$ is strict. 
\end{proof}
\subsection{Rees description of the spectral sequence}

$\hspace{3mm}$In the above setup, a filtered differential module $(\mr{M,F,d})$ gives rise to a differential Rees module $\left(\mr{R}_\mr{F}(\mr{M}),\mr{d}_\mr{R}\right)$. In this section, we explore the relation between this Rees module and the spectral sequence $(\mr{E}_r,\mr{d}_r)$ associated to the filtered differential module, where we omit the module from the notation of the spectral sequence if it's clear from the context.\\
$\hspace*{3mm}$ By \cref{basecase}, we know that the restriction of $\mr{R_\mr{F}(M)}$ to $\hbar=0$ is precisely the associated graded and we can express the spectral sequence as $$\mr{E}_1=\mr{H}\left(\mr{R_\mr{F}(M)}/\hbar\right)\Rightarrow \mr{R_\mr{F}(H(M))}/\hbar.$$ We can upgrade the former observation to cover all infinitesimal neighbourhoods of $\hbar=0$:

\begin{proposition}
For any $r\ge 1$, there is a natural isomorphism of $\mr{R_\mr{F}(A)}/\hbar^r$-modules $$\mr{R_\mr{F}(M)}/\hbar^r \cong \bigoplus_{p\in \Z} \mr{F}^p\mr{M}/\mr{F}^{p+r}\mr{M}\cdot\hbar^{-p}.$$
\end{proposition}
\begin{proof}
Very much the same as \cref{basecase}.
\end{proof}
In particular, the weight $p$ component of the truncation $\mr{R_\mr{F}(M)}/\hbar^r$ is 

$$\left(\mr{F}^p\mr{M}/\mr{F}^{p+r}\mr{M}\right)\cdot\hbar^{-p}$$ and the natural map

 $$\mr{F}^p\mr{M}/\mr{F}^{p+r}\mr{M}\to\mr{F}^{p-r+1}\mr{M}/\mr{F}^{p+1}\mr{M}$$ corresponds to

 $$\cdot\hbar^{r-1}: \left(\mr{F}^p\mr{M}/\mr{F}^{p+r}\mr{M}\right)\cdot\hbar^{-p} \to \left(\mr{F}^{p-r+1}\mr{M}/\mr{F}^{p+1}\mr{M}\right)\cdot\hbar^{-p+r-1}.$$
 
  Hence by \cref{ssiso}, we deduce:
 \2
 \begin{corollary}\label[corollary]{pages}
 Let $(\mr{M,F})$ be a filtered differential module with an associated spectral sequence $(\mr{E}_r,\mr{d}_r)$. For all $r\ge 1$ and all $p \in \Z$ there are natural isomorphisms
 $$\mr{E}_r^p \cong \mr{im}\left(\cdot \hbar^{r-1}: \mr{H}^p\left(\mr{R_\mr{F}(M)}/\hbar^r\right)\to \mr{H}^{p-r+1}\left(\mr{R_\mr{F}(M)}/\hbar^r\right)\right).$$
 \end{corollary}
 
$\hspace*{3mm}$ In the same fashion, following the discussion after \cref{ssiso}, we can describe the differentials $\mr{d}_r$ in terms of the Rees module and its differential. Consider the exact sequence of $k[\hbar]$-modules:

 $$0\xrightarrow{} \mr{R}_\mr{F}(\mr{M})/\hbar^r \xrightarrow{\cdot \hbar^r}  \mr{R}_\mr{F}(\mr{M})/\hbar^{2r} \xrightarrow{} \mr{R}_\mr{F}(\mr{M})/\hbar^r \xrightarrow{} 0,$$
 
where the action of $\hbar^r$ on the left is given by taking an element in $\mr{R}_\mr{F}(\mr{M})/\hbar^r$, lifting it to $\mr{R}_\mr{F}(\mr{M})$, multiplying by $\hbar^r$ and then taking the natural projection to the quotient $\mr{R}_\mr{F}(\mr{M})/\hbar^{2r}$. By naturality of the connecting morphism, we get a commutative square  \begin{equation}\label{diagram1}
\begin{tikzcd}
\mr{H}^p\left(\mr{R_\mr{F}(M)}/\hbar^r\right) \arrow[r, "\delta_{r,r}"] \arrow["\cdot\hbar^{r-1}",d]
& \mr{H}^{p+r}\left(\mr{R_\mr{F}(M)}/\hbar^r\right)\arrow[d,"\cdot\hbar^{r-1}"]   \\
\mr{H}^{p-r+1}\left(\mr{R_\mr{F}(M)}/\hbar^r\right)\arrow[r, "\delta_{r,r}"]
& \mr{H}^{p+1}\left(\mr{R_\mr{F}(M)}/\hbar^r\right).
\end{tikzcd}
\end{equation}
Hence the bottom horizontal map restricts to a map between the images of the vertical maps and we have the following:
 \2
 \begin{corollary}\label[corollary]{diffs}
 Let $(\mr{M,F})$ be a filtered differential module with an associated spectral sequence $(\mr{E}_r,\mr{d}_r)$. For all $r\ge 1$ and all $p \in \Z$ there are commutative diagrams:
 \[
\begin{tikzcd}
\mr{E}_r^p \arrow[r, "\mr{d}_{r}"] \arrow["\cong",d]
& \mr{E}_r^{p+r}\arrow[d,"\cong"]   \\
\mr{im}\left(\cdot \hbar^{r-1}|_{\mr{H}^p\left(\mr{R_\mr{F}(M)}/\hbar^r\right)}\right)\arrow[r, "\delta_{r,r}"]
& \mr{im}\left(\cdot \hbar^{r-1}|_{\mr{H}^{p+r}\left(\mr{R_\mr{F}(M)}/\hbar^r\right)}\right).
\end{tikzcd}
\]
In other words, under the isomorphisms of \cref{pages}, the differential of the spectral sequence $\mr{d}_r:\mr{E}^p_r\to \mr{E}^{p+r}_r$ corresponds to the restriction of the connecting morphism $\delta_{r,r}$ to the images of the vertical maps in the above commutative square (\ref{diagram1}).
 \end{corollary}
 \2
\begin{remark}We can give yet another description of the spectral sequence using the exact couple arising from the short exact sequence:
 
 $$0\to \mr{R_\mr{F}(M)}\xrightarrow{\cdot\hbar}\mr{R_\mr{F}(M)}\to \mr{R_\mr{F}(M)}/\hbar\to 0.$$
 \end{remark}

\section{Spectral sequences and deformation theory}
 \subsection{\texorpdfstring{$\mr{E}_1$}--degeneration and flatness}
 
 $\hspace*{3mm}$The description of the spectral sequence associated to a filtered differential module in terms of the Rees functor allows us to rephrase the classical $\mr{E}_1$-degeneration from \cref{e1} in more geometric fashion.
 \2
 \begin{proposition}\label[proposition]{e1flatness}
 Let $\mr{(M,F,d)}$ be a biregular filtered differential module. Its spectral sequence collapses on $\mr{E}_1$ if and only if the differential Rees module $\mr{\left(R_\mr{F}(M),d_R\right)}$ is flat, i.e. its cohomology is flat over $k[\hbar]$.
 \end{proposition}
 \begin{proof}
 By \cref{flatness} flatness is equivalent to the strictness of $\mr{d}$. For a biregular filtration, \cref{e1} implies the latter is, in turn, equivalent to $\mr{E}_1$-degeneration.
 \end{proof}
 
 $\hspace*{3mm}$Since $k[\hbar]$ is a PID, the structure theorem for finitely generated modules yields a decomposition 
 $$\mr{H}(\mr{R_\mr{F}(M)})\cong k[\hbar]^{\oplus r}\oplus \bigoplus_{i} k[\hbar]/\hbar^{r_i},$$ where for simplicity we drop the powers of $\hbar$ that track the filtration from the notation and are responsible for the exact locations of the cycles on the various pages. We can read off the entire spectral sequence from this -- the free part corresponds to permanent cycles that survive all pages. Each torsion module $k[\hbar]/\hbar^{r_i}$ corresponds to a pair of cycles on $\mr{E}_{r_i}$ and an identity differential between them. 
 \2
 \begin{remark}
 We can explain the strictness equation $$\mr{d}\left(\mr{F}^p\mr{M}\right)=\mr{d}(\mr{M})\cap \mr{F}^p\mr{M}$$ in terms of torsion. Namely, the failure of this equation to hold means that we have an $v \in \mr{F}^p\mr{M}$ such that $v=\mr{d}w$ but $w\in \mr{F}^{p-r}\mr{M}$ for some $r\ge 1$. On the Rees side, we get generators $e_v$ and $e_w$ such that $$\mr{d}_\mr{R}\left(e_w\right)=\hbar^re_v.$$ This creates torsion in cohomology and vice-versa. 
 \end{remark}
\subsection{Obstructions to deforming cohomology classes}

$\hspace*{3mm}$Let $(\mr{M,F},\mr{d})$ be a biregular filtered differential module. Applying the Rees construction, we get a perfect complex of graded $k[\hbar]$-modules $(\mr{R_\mr{F}(M)},\mr{d}_\mr{R})$.\\
$\hspace*{3mm}$ Recall that the spectral sequence associated to $\mr{M}$ can be expressed in terms of $\mr{R_\mr{F}(M)}$ instead as $$\mr{H(R_\mr{F}(M)}/\hbar)\Rightarrow \mr{R_\mr{F}(H(M))}/\hbar.$$ $\hspace*{3mm}$Consider the short exact sequence $$0\to \mr{R_\mr{F}(M)}\xrightarrow{\cdot\hbar}\mr{R_\mr{F}(M)}\to \mr{R_\mr{F}(M)}/\hbar\to 0.$$ It induces a long exact sequence in cohomology \begin{equation}\label{deltak}\cdots \longrightarrow \mr{H}^n(\mr{R}_\mr{F}(\mr{M))} \xrightarrow{\;\cdot \hbar\;} \mr{H}^n(\mr{R}_\mr{F}\mr{(M))} \xrightarrow{\;\pi_n\;} \mr{E}_1^n \xrightarrow{\;\delta_n\;} \mr{H}^{n+1}(\mr{R}_\mr{F}\mr{(M))} \xrightarrow{\;\cdot \hbar\;} \mr{H}^{n+1}(\mr{R}_\mr{F}\mr{(M))} \longrightarrow \cdots\end{equation} 
$\hspace*{3mm}$It follows that a class in $\mr{E}_1^n$ extends to $\mr{H}^n(\mr{R}_\mr{F}\mr{(M))}$ if and only if it is in $\mr{ker}(\delta_n)$. Hence all $\mr{E}_1^n$-classes extend if and only if $\delta_n=0$. The latter condition is equivalent to the map $$\mr{H}^n(\mr{R}_\mr{F}\mr{(M)}) \xrightarrow{\;\cdot \hbar\;} \mr{H}^n(\mr{R}_\mr{F}\mr{(M))}$$ having trivial kernel, i.e. the $k[\hbar]$-modules $\mr{H}^n(\mr{R}_\mr{F}\mr{(M))}$ are flat (and thus free). Equivalently, we can split the above long exact sequence into short exact sequences 
\begin{equation}\label{torss}0 \longrightarrow \mr{H}^n(\mr{R}_\mr{F}\mr{(M)})/\hbar \longrightarrow \mr{E}_1^n \longrightarrow \mr{Tor}^{k[\hbar]}_1\left(\mr{H}^{n+1}(\mr{R}_\mr{F}(\mr{M)}),k\right)\longrightarrow 0.\end{equation}
$\hspace*{3mm}$Thus $\mr{E}_1^n$-classes which do not extend correspond precisely to torsion elements in $\mr{H}^{n+1}(\mr{R}_\mr{F}\mr{(M)})$, hance by \cref{e1flatness} we obtain the following result.
\2
\begin{proposition}\label[proposition]{lifts}
Let $(\mr{M,F},\mr{d})$ be a biregular filtered differential module. Its spectral sequence degenerates at $\mr{E}_1$ if and only if any $\mr{E}_1$-class extends to a class in $\mr{H}\left(\mr{R_\mr{F}(M)}\right)$.
\end{proposition}
\2
\begin{remark}
In the discussion above, if we work with a filtered complex $\mr{C}$ instead, then both the spectral sequence and the Rees complex are naturally bigraded. The terms of the Rees complex are graded by the filtration: $$\mr{R}^n_\mr{F}\mr{(C)}=\bigoplus_{p\in \Z}\mr{R}^{p,n-p}_\mr{F}\mr{(C)},$$where $$\mr{R}^{p,n-p}_\mr{F}\mr{(C)} = \mr{F}^p\mr{C}^n\hbar^{-p}.$$ 
$\hspace*{3mm}$As the differential of $\mr{C}$ is compatible with the filtration, the Rees differential preserves the weight and hence the Rees cohomology is bigraded: $$\mr{H}^n(\mr{R}_\mr{F}(\mr{C)})=\bigoplus_{p+q=n}\mr{H}^{p,q}(\mr{R}_\mr{F}\mr{(C)}),$$ where $$\mr{H}^{p,q}(\mr{R}_\mr{F}\mr{(C)})=\mr{H}^{p+q}(\mr{R}^{p,\bullet-p}_\mr{F}\mr{(C)}).$$ 
$\hspace*{3mm}$For $x\in \mr{F}^p\mr{C}^n$, we have $$\hbar\cdot x\hbar^{-p}=x\hbar^{-(p-1)} \in \mr{F}^{p-1}\mr{C}^n\hbar^{-(p-1)},$$ so multiplication by $\hbar$ has bidegree $(-1,1)$, i.e. we have $$\hbar\cdot :\mr{H}^{p,q}(\mr{R}_\mr{F}\mr{(C)}) \to \mr{H}^{p-1,q+1}(\mr{R}_\mr{F}\mr{(C)}).$$ 
The short exact sequence $$0\to \mr{R}^{p+1,\bullet}_\mr{F}\mr{(C)}\xrightarrow{\cdot\hbar}\mr{R}^{p,\bullet}_\mr{F}\mr{(C)}\to \mr{R}^{p,\bullet}_\mr{F}\mr{(C)}/\hbar\to 0$$ yields as above $$0 \longrightarrow \mr{H}^{p,q}(\mr{R_\mr{F}(C)})/\hbar \longrightarrow \mr{E}_1^{p,q} \longrightarrow \mr{Tor}^{k[\hbar]}_1\left(\mr{H}^{p,q+1}(\mr{R_\mr{F}(C)}),k\right)\longrightarrow 0.$$ 
\end{remark}
\subsection{Finite order deformations}
$\hspace*{3mm}$Extending all the way to $\mr{H}^n(\mr{R_\mr{F}(C)})$ at once might be quite difficult in practice. Instead, we could work with the truncated/pro Rees complexes, extending inductively over the infinitesimal neighbourhoods of the origin. In addition, the biregularity condition is unnecessarily strong as we shall see below.\\
$\hspace*{3mm}$Recall that an element in $\mr{R}_\mr{F,N}\mr{(C)}=\mr{R_\mr{F}(C)}/\hbar^{N+1}$ of total degree $n$, lifting $x_0 \in \mr{E}_1^{p,n-p}$, has the form $$\mathbf{x}^{(N)}_\hbar=x_0\hbar^{-p}+x_1\hbar^{-p}+\cdots+x_N\hbar^{-p},$$ where $x_r\in \mr{F}^{p+r}\mr{C}^n$. This is a cocycle if and only if $$\mr{d_R}\left(\mathbf{x}^{(N)}_\hbar\right)=\left(\mr{d}x_0+\mr{d}x_1+\cdots+\mr{d}x_N\right)\hbar^{-p}\in \hbar^{N+1}\mr{R}^{n+1}_\mr{F}\mr{(C)}.$$ 
$\hspace*{3mm}$In terms of the filtration, we get the condition $$\mr{d}\left(x_0+x_1+\cdots+x_N\right)\in \mr{F}^{p+N+1}\mr{C}^{n+1}.$$
\begin{theorem}\label{infinitesimal}
Let $(\mr{C,F},\mr{d})$ be a filtered complex with associated spectral sequence $\left(\mr{E}_r,\mr{d}_r\right)$. A class $x_0$ on the central fibre of the Rees complex, that is $x_0\in \mr{E}_1$, extends to the $N^{\text{th}}$-infinitesimal neighbourhood if and only if $$\mr{d}_1\left([x_0]\right)=0,\;\; \mr{d}_2\left([x_0]\right)=0,\cdots, \mr{d}_N\left([x_0]\right)=0.$$
\end{theorem}
\begin{proof}
Take a class $[x_0] \in \mr{E}^{p,n-p}_1$ and $\mr{d}(x_0)\in \mr{F}^{p+1}\mr{C}^n$. Then $\mr{d}_1([x_0])$ is represented by $$\mr{d}(x_0) \mod \mr{F}^{p+2}\mr{C}^{n+1}.$$ So $\mr{d}_1([x_0])=0$ if and only if there exists $x^{(1)} \equiv x_0\mod \mr{F}^{p+1}$ and $\mr{d}\left(x^{(1)}\right) \in \mr{F}^{p+2}\mr{C}^{n+1}$. Hence, after $N$ steps, we have $$\mr{d}_1(x_0)=0,\cdots,\mr{d}_N(x_0)=0$$ if and only if there exists $$x^{(N)} \equiv x_0\mod \mr{F}^{p+1}\text{ with } \mr{d}\left(x^{(N)}\right) \in \mr{F}^{p+N+1}\mr{C}^{n+1}.$$ 
$\hspace*{3mm}$When lifting from $r$ to $r+1$, we may have to change the lower order corrections. Having chosen these corrections by induction, we may write $x^{(N)} = x_0+x_1+\cdots +x_N$ with $x_r \in \mr{F}^{p+r}\mr{C}^{n}$, the corresponding Rees element is $$\mathbf{x}^{(N)}_\hbar\coloneqq x^{(N)}\hbar^{-p}= \sum_{r=0}^N \left(x_r\cdot\hbar^{-(p+r)}\right)\hbar^{r}.$$ 
$\hspace*{3mm}$The condition $\mr{d}\left(x^{(N)}\right) \in \mr{F}^{p+N+1}\mr{C}^{n+1}$ translates to $\mr{d}_\mr{R}\left(\mathbf{x}^{(N)}_\hbar\right)\in \hbar^{N+1}\mr{R}^{n+1}_\mr{F}(\mr{C})$.
\end{proof}
\2
\begin{corollary}
Let $(\mr{C,F},\mr{d})$ be a filtered complex with associated spectral sequence $\left(\mr{E}_r,\mr{d}_r\right)$. A class $x_0$ on the central fibre of the Rees complex, that is $x_0\in \mr{E}_1$, extends to the pro Rees complex $\mathbf{R}_\mr{F}(\mr{C})$ if and only if $\mr{d}_r\left([x_0]\right)=0$ for all $r\ge 1$.
\end{corollary}
\begin{proof}
This follows immediately from the definitions and \cref{infinitesimal}.
\end{proof}
$\hspace*{3mm}$In addition to allowing us to use induction in practice, this shows that the biregularity condition on the filtration can be relaxed. For example, flatness implies $\mr{E}_1$-degeneration under no extra conditions.
\2
\begin{corollary}\label[corollary]{flate1}
Let $(\mr{C,F},\mr{d})$ be a filtered complex. Assume that $\mr{H}(\mr{R}_\mr{F}(\mr{C}))$ is flat over $k[\hbar]$. Then the spectral sequence degenerates at $\mr{E}_1$.
\end{corollary}
\begin{proof}
We have a short exact sequence $$0 \longrightarrow \mr{H}^{n}\left(\mr{R}_\mr{F}(\mr{C)}\right)/\hbar \longrightarrow \mr{E}_1^n \longrightarrow \mr{Tor}^{k[\hbar]}_1\left(\mr{H}^{n+1}\left(\mr{R}_\mr{F}(\mr{C)}\right),k\right)\longrightarrow 0$$ the last term $$\mr{Tor}^{k[\hbar]}_1\left(\mr{H}^{n+1}(\mr{R}_\mr{F}(\mr{C})),k\right)=0$$ by flatness, hence every $\mr{E}_1$-class comes from $\mr{H}(\mr{R}_\mr{F}(\mr{C)})/\hbar$. In particular, it extends to a class in $\mr{H}(\mr{R}_\mr{F}(\mr{C)})$, thus it lifts to the $N^{\text{th}}$-infinitesimal neighbourhood for any $\mr{N}\ge 1$, so the result follows from \cref{infinitesimal}.
\end{proof}
$\hspace*{3mm}$For the converse, without biregularity, we need an extra condition:
\2
\begin{corollary}\label[corollary]{separatede1}
Let $(\mr{C,F},\mr{d})$ be a filtered complex and suppose its spectral sequence collapses at $\mr{E}_1$. Assume that the Rees cohomology $\mr{H}(\mr{R}_\mr{F}(\mr{C)})$ is $\hbar$-adically separated: $$\bigcap_{N\ge 1} \hbar^{N}\mr{H}(\mr{R}_\mr{F}(\mr{C)})=0.$$ Then $\mr{H}(\mr{R}_\mr{F}(\mr{C)})$ is flat.
\end{corollary}
\begin{proof}
Indeed, every $x_0 \in \mr{E}_1^n$-class lifts to the $N^{\text{th}}$-infinitesimal neighbourhood for any $N\ge 1$. Choose a representative $x_N \in \mr{R}^n_\mr{F}\mr{(C)}$ with $$\mr{d}\left(x_N\right)=\hbar^{N}y_N.$$ 
Hence $$\delta_n\left(x_0\right)=\hbar^N\left[y_N\right]\in \mr{H}^n(\mr{R}_\mr{F}(\mr{C)})$$ and consequently $$\delta_n\left(x_0\right)\in \hbar^N\mr{H}^n_\mr{R}(\mr{C,F})$$ for all $N\ge 1$. Since this is true for every class $x_0$, we have $$\mr{im}\left(\delta_n\right)\subseteq   \bigcap_{N\ge 1} \hbar^{N}\mr{H}^n(\mr{R}_\mr{F}(\mr{C)})=0.$$ By \cref{deltak}, we have $$\mr{Tor}^{k[\hbar]}_1\left(\mr{H}^{n+1}(\mr{R}_\mr{F}\mr{(C))},k\right)=\mr{im}\left(\delta_n\right),$$ which finishes the proof.
\end{proof}
\subsection{Formal deformations}
$\hspace*{3mm}$In the biregular case, there are no compatibility obstructions.
\2
\begin{proposition}
The spectral sequence $\left(\mr{E}_r,\mr{d}_r\right)$ of a biregular $(\mr{C,F},\mr{d})$ degenerates on $\mr{E}_1$ if and only if all classes on the central fibre $\mr{E}_1$ extend to $\mr{H}^n\left(\widehat{\mr{R}}_\mr{F}(\mr{C)}\right)$ over the formal neighbourhood of the origin $\widehat{\mathbf{A}^1}_0$.
\end{proposition}
\begin{proof}
If $\mr{E}_1=\mr{E}_\infty$, by \cref{lifts} every $\mr{E}_1$ class extends to a class in $\mr{H}(\mr{R}_\mr{F}(\mr{C)})$. Completing the latter gives a formal extension.\\
$\hspace*{3mm}$Conversely, suppose any $\mr{E}_1$ class extends to $\mr{H}^n\left(\widehat{\mr{R}}_\mr{F}(\mr{C)}\right)$, so the composition $$\mr{H}^n\left(\widehat{\mr{R}}_\mr{F}(\mr{C)}\right) \to \widehat{\mr{H}}^n(\mr{R}_\mr{F}(\mr{C)}) \to \mr{H}^n(\mr{R}_\mr{F}(\mr{C)})/\hbar \to \mr{E}^n_1$$ is surjective, in particular the last map $$\mr{H}^n(\mr{R}_\mr{F}(\mr{C})/\hbar \to \mr{E}^n_1$$ must be surjective, hence by \cref{torss} the Rees cohomology is flat and by \cref{flate1} the spectral sequence degenerates at $\mr{E}_1$.
\end{proof}
$\hspace*{3mm}$In general, we have:
\2
\begin{theorem}\label{formale1}
Let $(\mr{C,F},\mr{d})$ be a filtered complex. Assume that $\mr{H}\left(\widehat{\mr{R}}_\mr{F}(\mr{C})\right)$ is flat over $k\brak$; equivalently, every $\mr{E}_1$-class has a formal lift. Then the spectral sequence degenerates at $\mr{E}_1$.\\
Conversely, if the spectral sequence degenerates at $\mr{E}_1$ and $\mr{H}\left(\widehat{\mr{R}}_\mr{F}(\mr{C})\right)$ satisfies the separateness condition: $$\bigcap_{N\ge 1} \hbar^{N}\mr{H}\left(\widehat{\mr{R}}_\mr{F}(\mr{C)}\right)=0,$$ then $\mr{H}\left(\widehat{\mr{R}}_\mr{F}(\mr{C})\right)$ is flat and every $\mr{E}_1$-class has a formal lift.
\end{theorem}
\2
\begin{remark}\label[remark]{relativermk}
The results extend naturally to the relative setting over $\mr{R_F(A)}$, without any changes, except that flatness has to be replaced with vanishing of $\hbar$-torsion which we may call $\hbar$-flat/flat along the Rees parameter.	
\end{remark}
\subsection{Sheaves over affine schemes}
$\hspace*{3mm}$Let $\mathbf{A}^1=\mr{Spec}(k[\hbar])$ be the affine line which we equip with the standard action of the multiplicative group $\mathbf{G}_m$. 
\2 
\begin{definition}
 Let $(\mr{A},\mr{F})$ be a filtered algebra. The Rees scheme of $(\mr{A,F})$ is the $\mathbf{A}^1$-scheme $\mr{Spec}(\mr{R}_\mr{F}(\mr{A})).$
\end{definition}
$\hspace*{3mm}$The grading of $\mr{R}_\mr{F}(\mr{A})$ corresponds to a $\mathbf{G}_m$-action on $\mr{Spec}(\mr{R}_\mr{F}(\mr{A}))$.
\2
\begin{definition}
Let $(\mr{M,F})$ be a filtered module over a filtered algebra $(\mr{A},\mr{F})$. The Rees bundle of $(\mr{M,F})$ is the $\mathbf{G}_m$-equivariant sheaf $\mathscr{R}_\mr{F}(\mr{M})$ over the Rees scheme $\mr{Spec}(\mr{R}_\mr{F}(\mr{A}))$ associated to the graded $\mr{R}_\mr{F}(\mr{A})$-module $\mr{R}_\mr{F}(\mr{M})$.
\end{definition}

\2
\begin{proposition}
 Let $(\mr{A},\mr{F})$ be a filtered algebra. The Rees sheaf construction induces an exact functor $$\mathscr{R}: \mr{Fil}(\mr{Mod}(\mr{R}_\mr{F}(\mr{A}))\to \mr{QCoh}^{\mathbf{G}_m}(\mr{Spec}(\mr{R}_\mr{F}(\mr{A}))$$ which is compatible with direct sums and tensor products.
\end{proposition}
\begin{proof}
Indeed, $\mathscr{R}$ is the composition $\widetilde{(-)} \circ \mr{R}$, so the claim follows.
\end{proof}
\2
\begin{proposition}
 Let $j : \mathbf{G}_m \xhookrightarrow{} \mathbf{A}^1$ be the inclusion. Let $\mr{(A,F)}$ be a filtered algebra and let $\mr{(M,F)}$ be a filtered module over it. There are natural isomorphisms:
 \begin{enumerate}
  \item $\mathscr{R}_\mr{F}(\mr{M})_{1} \simeq \mr{M}$ as $\mr{A}$-modules.
    \item $j^{*}\mathscr{R}_\mr{F}(\mr{M}) \simeq \mathscr{O}_{\mathbf{G}_m}\otimes \mr{M}$ as $\mathscr{O}_{\mathbf{G}_m}\otimes \mr{A}$-modules.
  \item $\mathscr{R}_\mr{F}(\mr{M})_{0} \simeq \mr{Gr}_\mr{F}\mr{M}$ as graded $\mr{Gr_F(A)}$-modules.
 \end{enumerate}
\end{proposition}
\begin{proof}
It follows from \cref{basecase} since the functor $\widetilde{(-)}$ is compatible with localisation and pullback.
\end{proof}
$\hspace*{3mm}$Given a $\mathbf{G}_m$-equivariant vector bundle $\mathscr{V}$ on $\mathbf{A}^1$, let $\mr{V} = \mathscr{V}_1$. On the other hand, $\mathscr{V}$ corresponds to a free $k[\hbar]$-module via global sections $\Gamma(\mathbf{A}^1,\mathscr{V})$. The trivialisation over $\mathbf{G}_m \subset \mathbf{A}^1$ gives an isomorphism $$\Gamma(\mathbf{A}^1,\mathscr{V})\otimes_{k[\hbar]}k[\hbar,\hbar^{-1}] \simeq \mr{V}\otimes k[\hbar,\hbar^{-1}].$$ Hence, $\Gamma(\mathbf{A}^1,\mathscr{V}) \subset \mr{V}\otimes k[\hbar,\hbar^{-1}]$ and we get a filtration on $\mr{V}$ by setting $$\mr{F}^p\mr{V} = \{s \in \mr{V}, \hbar^{-p}s \in \Gamma(\mathbf{A}^1,\mathscr{V})\}.$$
This construction is (quasi-)inverse to the Rees bundle construction, and we get:
\2
\begin{proposition}
 There is an equivalence between filtered vector spaces and $\mathbf{G}_m$-equivariant vector bundles on $\mathbf{A}^1$.
\end{proposition}

$\hspace*{3mm}$We can rephrase the results obtained in the previous sections in a geometric/sheaf-theoretic language. Recall that given a filtered complex $\mr{(C,F)}$, we have a sheafified Rees complex $$\mathscr{R}_\mr{F}(\mr{C})\in \mr{Perf}(\mathbf{A}^1).$$ 
$\hspace*{3mm}$Let $$\mr{D}_N=\mr{Spec}\left(k[\hbar]/\hbar^{N+1}\right)$$ be the $N^{\text{th}}$ infinitesimal neighbourhood of the origin.
\2
\begin{proposition}\label[proposition]{geometric}
Let $(\mr{C,F},\mr{d})$ be a biregular filtered complex. The following are equivalent:
\begin{enumerate}
\item The spectral sequence $(\mr{E}_r,\mr{d}_r)$ degenerates at $\mr{E}_1$.
\item The sheaves $\mr{H}^n(\mathscr{R}_\mr{F}(\mr{C})$ are flat $\mathscr{O}_{\mathbf{A}^1}$-modules.
\item For all $n\in \Z$, any cohomology class on the central fibre $\mr{H}^n(\mathscr{R}_\mr{F}(\mr{C})_0)$ extends to a global section of $\mr{H}^n(\mathscr{R}_\mr{F}(\mr{C)})$ over $\mathbf{A}^1$.
\item For all $n\in \Z$, any cohomology class on the central fibre $\mr{H}^n(\mathscr{R}_\mr{F}(\mr{C})_0)$ extends to a global section of $\mr{H}^n\left(\mathscr{R}_\mr{F}(\mr{C)}|_{\mr{D}_{N}}\right)$ over the infinitesimal neighbourhood $\mr{D}_N$ for all $N\ge 1$.
\item For all $n\in \Z$, any cohomology class on the central fibre $\mr{H}^n\left(\mathscr{R}_\mr{F}(\mr{C})_0\right)$ extends to a section of $\mr{H}^n\left(\widehat{\mathscr{R}}_\mr{F}(\mr{C)}\right)$ over the formal neighbourhood $\widehat{\mathbf{A}^1}_0$.
\end{enumerate}
\end{proposition}
\2
\begin{remark}
In \cref{geometric}, flatness implies $\mr{E}_1$-degeneration for any filtered complex and we can relax the biregularity condition on the filtration for the opposite direction by introducing the separatedness condition on the cohomology as above. 
\end{remark}
\2
\begin{remark}
Again, everything holds true in the relative setting over $\mathcal{X} = \mr{Spec}\left(\mr{R_F(A)}\right)$ except that flatness has to be replaced with the $\mr{Tor}$ vanishing condition $$\mr{Tor}_1^{\mathscr{O}_\mathcal{X}}\left(\mathscr{R}_\mr{F}(\mr{M}),\mathscr{O}_{\mathcal{X}_0}\right)=0$$ which corresponds to the vanishing of $\hbar$-torsion mentioned in \cref{relativermk}.
\end{remark}
$\hspace*{3mm}$We conclude with the following formality result for the Rees complex.
\2
\begin{corollary}
Let $(\mr{C,F},\mr{d})$ be a filtered complex and suppose that $(\mr{E}_r,\mr{d}_r)$ degenerates at $\mr{E}_1$.
\begin{enumerate}
\item There exists a quasi-isomorphism $$\mathscr{R}_\mr{F}(\mr{C}) \simeq \bigoplus_n \mr{H}^n(\mathscr{R}_\mr{F}(\mr{C)})[-n].$$ In other words, the Rees complex is formal in $\mr{Perf}\left(\mathbf{A}^1\right)$. 
\item The completed Rees complex $\widehat{\mathscr{R}}_\mr{F}(\mr{C})$ is formal in $\mr{Perf}\left(\widehat{\mathbf{A}^1}_0\right)$.
\end{enumerate}
\end{corollary}
\begin{proof}
We prove the first claim, the same proof applies to the second. Since the spectral sequence collapses at $\mr{E}_1$, it follows by \cref{geometric} that the cohomology is free and hence projective. Since the higher $\mr{Ext}^{\ge 2}$ vanish, we obtain the claim.
\end{proof}
\subsection{Multiplicative structures}
\label{mult}
$\hspace*{3mm}$Consider a filtered dg algebra $(\mr{A,F,d})$. Since the Rees functor is monoidal, we obtain a filtered dg algebra $\mr{R}_\mr{F}(\mr{A})$ over $k[\hbar]$. Using the multiplication, it is enough to lift a set of generators.
\2
\begin{proposition}
Suppose $(\mr{A,F,d})$ is a filtered dg algebra and assume that the algebra $\mr{H}(\mr{Gr}_\mr{F}(\mr{A}))$ is generated by classes that admit lifts to $\mr{H}(\mr{R}_\mr{F}(\mr{A}))$. Then $\mr{H}(\mr{R}_\mr{F}(\mr{A}))$ is flat over $k[\hbar]$ and hence the spectral sequence associated to $\mr{(A,F,d)}$ collapses at $\mr{E}_1$. 
\end{proposition}
\begin{proof}
Since $\mr{H}(\mr{Gr}_\mr{F}(\mr{A}))$ is generated by classes that extend to $\mr{H}(\mr{R}_\mr{F}(\mr{A}))$, we see that the natural map $$\mr{H}(\mr{R}_\mr{F}(\mr{A})) \to \mr{H}(\mr{Gr}_\mr{F}(\mr{A}))$$ is surjective, hence $\mr{H}(\mr{R}_\mr{F}(\mr{A}))$ has no $\hbar$-torsion and is therefore flat. The result follows from \cref{flate1}
\end{proof}
$\hspace*{3mm}$Consider a compact K\"ahler Lagrangian $\mr{L}$ in a holomorphic symplectic variety $\mr{X}$. We apply this observation to the local-to-global $\mr{Ext}$ spectral sequence of a line bundle $\mathscr{L}$ on $\mr{L}$. 
\2
\begin{corollary}\label[corollary]{e2lag}
Let $\mr{X}$ be a holomorphic symplectic variety and let $i:\mr{L}\xhookrightarrow{} \mr{X}$ be a compact K\"ahler Lagrangian whose Hodge algebra $$\bigoplus_{p,q}\mr{H}^p(\Omega_\mr{L}^q)$$ is generated by classes of holomorphic degree $q\le 1$. Suppose $\mathscr{L}$ is a line bundle on $\mr{L}$. Then the $\mr{Ext}$ spectral sequence $$\mr{H}^p(\Omega_\mr{L}^q) \Rightarrow \mr{Ext}^{p+q}(i_*\mathscr{L},i_*\mathscr{L})$$ collapses at $\mr{E}_2$ if and only if $\mr{d}_2^{p,1}=0$ for all $p\ge 0$.
\end{corollary}
\2
\begin{remark}
This is similar to some of the results in \cite{MladenovDegeneration}. It requires more than just $\mr{d}_2(\omega)=0$, but doesn't require that $\mathscr{L}$ extends to first order in every direction. 
\end{remark}
$\hspace*{3mm}$In particular, \cref{e2lag} works out very well in the case of an abelian Lagrangian. Since the Hodge cohomology is generated by $\mr{H}^{0,1}(\mr{L})\oplus\mr{H}^{1,0}(\mr{L})$, we only need to check the vanishing of a single differential:
\2
\begin{corollary}\label[corollary]{abeliancase}
Let $\mr{X}$ be a holomorphic symplectic variety and let $i:\mr{L}\xhookrightarrow{} \mr{X}$ be an abelian Lagrangian. Let $\mathscr{L}$ be a line bundle on $\mr{L}$. Then the $\mr{Ext}$ spectral sequence $$\mr{H}^p(\Omega_\mr{L}^q) \Rightarrow \mr{Ext}^{p+q}(i_*\mathscr{L},i_*\mathscr{L})$$ collapses at $\mr{E}_2$ if and only if $\mr{d}_2^{0,1}=0$.
\end{corollary}
\2
\begin{remark}
The vanishing in \cref{abeliancase} hold for example for $\mathscr{L}=\mathscr{O}_\mr{L}$.
\end{remark}
\subsection{Hodge-to-de Rham for K\"{a}hler manifolds and beyond}

$\hspace*{3mm}$Let $\mr{X}$ be a compact K\"{a}hler manifold of dimension $n$. Denote the sheaf of smooth differential $k$-forms on $\mr{X}$ by $\mathcal{A}^k_\mr{X}$. Let $$\mr{A}^k_{\mr{X}}=\Gamma\left(\mr{X},\mathcal{A}^k_\mr{X}\right)$$ be the space of global sections. The de Rham differential decomposes $\mr{d}_\mr{dR}=\partial +\bar{\partial}$ and we get a decomposition $$\mathcal{A}^k_\mr{X}=\bigoplus_{p+q=k}\mathcal{A}^{p,q}_\mr{X}$$ into $(p,q)$-types and similarly for $\mr{A}^k_\mr{X}$. We consider the Hodge filtration $$\mr{F}^p \mr{A}^k_\mr{X} = \bigoplus_{r \ge p} \Gamma\left(\mr{X}, \mathcal{A}^{r, k-r}_\mr{X}\right).$$
$\hspace*{3mm}$The Rees complex $\mr{R}_\mr{X} \coloneqq \mr{R}(\mr{A}^k_{\mr{X}},\mr{F})$ is a complex of graded $\C[\hbar]$-modules
$$\mr{R}_\mr{X}^k = \bigoplus_{p=0}^k \mr{F}^p\left(\mr{A}^k_{\mr{X}}\right) \hbar^{-p}.$$
$\hspace*{3mm}$From our previous discussion, its differential is $$\mr{d}_{\mr{R}_\mr{X}}=\bar{\partial}+\hbar\partial,$$
 and the cohomology of the central fibre recovers the Dolbeault cohomology of the manifold
$$\mr{H}(\mr{R}_\mr{X} / \hbar \mr{R}_\mr{X}) \cong \bigoplus_{p,q} \mr{H}^q(\mr{X}, \Omega^{p}_\mr{X}),$$ where $\Omega^{p}_\mr{X}$ is the sheaf of holomorphic $p$-forms, while the generic fibre gives de Rham cohomology of $\mr{X}$:
$$\mr{H}(\mr{R}_\mr{X}\otimes_{\C\brak}\C\cbrak) \cong \mr{H}_\mr{dR}(\mr{X},\C\cbrak).$$
$\hspace*{3mm}$According to \cref{e1flatness}, the collapse of the Hodge-to-de Rham spectral sequence on $\mr{E}_1$ is equivalent to the flatness of $\mr{H}(\mr{R}_\mr{X})$ over $\C[\hbar]$. This is in turn equivalent to the surjectivity of the projection $$\pi : \mr{H}(\mr{R}_\mr{X})\to \mr{H}(\mr{R}_\mr{X} / \hbar \mr{R}_\mr{X}).$$
$\hspace*{3mm}$Let $$\alpha \in \Gamma\left(\mr{X},\mathcal{A}^{p,q}_\mr{X}\right)$$ be a Dolbeault cocycle, i.e. $\bar{\partial}\alpha=0$. We would like to find a Rees cocycle $$\alpha_\hbar = \alpha\hbar^{-p} + \alpha_1\hbar^{-p} + \alpha_2\hbar^{-p} +\cdots,$$ where $\alpha_r \in \Gamma\left(\mr{X},\mathcal{A}^{p+r,q-r}_\mr{X}\right)$. In terms of homogeneous Rees elements, we have $$\alpha_\hbar = \widetilde{\alpha} + \widetilde{\alpha_1}\hbar + \widetilde{\alpha_2}\hbar^2 +\cdots,$$ where $\widetilde{\alpha}_r = \alpha_r\hbar^{-(p+r)}\cdot \hbar^r$, as $\alpha_r$ is in the $(p+r)^\text{th}$ step of the filtration. Expanding the equation $$(\bar{\partial}+\hbar\partial)\alpha_\hbar=0$$
leads to a system of PDEs. The free term gives $\bar{\partial}\alpha=0$ which is by assumption. The next $\hbar$ term gives $$\partial\alpha+\bar{\partial}\alpha_1=0.$$ Similarly, $\hbar^{2}$ leads to $$\partial\alpha_1+\bar{\partial}\alpha_2=0$$ and so on. To solve the first equation, we want to show that $\partial\alpha$ is $\bar{\partial}$-exact. It is certainly $\partial$-closed and since $\partial$ and $\bar{\partial}$ anti-commute, we see it is $\bar{\partial}$-closed since $\alpha$ is. Hence by the $\partial\bar{\partial}$-lemma implies that $$\partial \alpha$$ is $\partial\bar{\partial}$-exact. Hence, there is a smooth $(p,q-1)$-form $\beta$ such that $$\partial\alpha= \partial\bar{\partial}\beta.$$
$\hspace*{3mm}$Hence, we let $\alpha_1 = \partial\beta$. Since this is $\partial$-closed, we may choose $\alpha_2=\alpha_3=\cdots=0$ and obtain the lift $$\alpha_\hbar = \alpha\hbar^{-p}+\partial\beta\hbar^{-p}= \widetilde{\alpha}+\widetilde{\partial\beta}\hbar.$$
$\hspace*{3mm}$ The first example for which collapse at $\mr{E}_1$ fails is given by the Iwasawa manifold $\mr{G}/\Gamma$, where $\mr{G}$ is the group of complex upper-triangular $3\times 3$ matrices and $\Gamma$ is the subgroup generated by matrices with non-diagonal entries in the Gaussian integers $\Z[i]$. There are $3$ left-invariant holomorphic $1$-forms $\omega_1,\omega_2,\omega_3$ such that $$\mr{d}_\mr{dR}\left(\omega_1\right)=0, \;\;\; \mr{d}_\mr{dR}\left(\omega_2\right)=0, \;\;\; \mr{d}_\mr{dR}\left(\omega_3\right)=\omega_1\wedge \omega_2.$$
At the Rees level we have $$\left[\omega_1\wedge\omega_2\cdot \hbar^{-2}\right]\not=0 \text{ and }
\mr{d}_\mr{R_X}\left(\omega_3\hbar^{-1}\right)= \bar{\partial}\omega_3\hbar^{-1}+\hbar\partial\omega_3\hbar^{-1} = \hbar\left(\omega_1\wedge\omega_2\cdot\hbar^{-2}\right).$$ This means $$\hbar\cdot\left[\omega_1\wedge\omega_2\cdot \hbar^{-2}\right]=0$$ and corresponds to torsion in the cohomology which is precisely the non-triviality of the spectral sequence. In fact, the cyclic module $k[\hbar]/\hbar$ associated with this torsion arises for a length $2$ horizontal zigzag given by $\omega_3$ and $\omega_1\wedge\omega_2$.

\section{Hard Lefschetz and flatness}
\label{hlrf}

$\hspace*{3mm}$One of the main advantages of the Rees point of view is that degeneration of
a spectral sequence becomes vanishing of torsion in cohomology which can be studied without computing the higher differentials. In
geometric situations, the main tools are Hard Lefschetz and Serre duality.\\ 
$\hspace*{3mm}$We formulate the argument first over a general principal ideal domain,
and conclude with an application to the Rees complex.

\subsection{Framework}

Let
$R=k[\hbar]$
or, more generally, let $R$ be a PID. Let $\mr{K}\in \mathrm{Perf}(R)$. We assume that its cohomology modules
$\mr{H}^i(\mr{K})$ are finitely generated over $R$, and write
\[
\mr{T}^i:=\mr{H}^i(\mr{K})_{\mathrm{tor}}
\]
for their torsion submodules. We shall use two structures on $\mr{K}$.

\begin{enumerate}
\item Serre duality of dimension $2n$,
\begin{equation}
\label{self-duality-HL}
\mr{K}\simeq \mr{R}\operatorname{Hom}_R(\mr{K},R)[-2n].
\end{equation}
\item A degree two operator
\[
\mr{L}:\mr{K}\longrightarrow \mr{K}[2]
\]
whose action on cohomology satisfies the Hard Lefschetz property: for every
$i\leq n$,
\begin{equation}
\label{integral-HL}
\mr{L}^{\,n-i}:\mr{H}^i(\mr{K})\xrightarrow{\ \sim\ }\mr{H}^{2n-i}(\mr{K}).
\end{equation}
\end{enumerate}

\begin{remark}
The argument below uses only the Hard Lefschetz isomorphisms on cohomology and
the self-duality \eqref{self-duality-HL}. It does not require a
pre-existing $\mathfrak{sl}_2$-action, nor does it require an adjointness
relation between $\mr{L}$ and the duality pairing. In applications an
$\mathfrak{sl}_2$-structure is often the natural source of
\eqref{integral-HL}, but it is stronger than what is needed for the
torsion argument.
\end{remark}

$\hspace*{3mm}$Because $R$ has global dimension one, the hyper-Ext spectral sequence for
the right-hand side of \eqref{self-duality-HL} degenerates to short exact
sequences
\begin{equation*}
\label{UCT-HL}
0\longrightarrow
\operatorname{Ext}^1_R\bigl(\mr{H}^{2n+1-i}(\mr{K}),R\bigr)
\longrightarrow \mr{H}^i(\mr{K})
\longrightarrow
\operatorname{Hom}_R\bigl(\mr{H}^{2n-i}(\mr{K}),R\bigr)
\longrightarrow0.
\end{equation*}
The last term is torsion-free. Hence the torsion submodule of $\mr{H}^i(\mr{K})$ is precisely
\begin{equation*}
\label{torsion-duality-identification}
\mr{T}^i\simeq
\operatorname{Ext}^1_R\bigl(\mr{H}^{2n+1-i}(\mr{K}),R\bigr).
\end{equation*}
Indeed, we have:
\2
\begin{lemma}
For a finitely generated module $\mr M$ over a principal ideal domain,
\[
\operatorname{Ext}^1_R(\mr{M},R)=0
\quad\text{if and only if}\quad
\mr{M}\text{ is torsion-free}.
\]
\end{lemma}
\begin{proof}
By the structure theorem for finitely generated modules over a PID, we can write
\[
\mr{M}\simeq R^{\oplus r}\oplus\bigoplus_a R/(m_a),
\]
taking $\mr{Exts}$, we obtain
\[
\operatorname{Ext}^1_R(\mr{M},R)\simeq\bigoplus_a R/(m_a).
\] which proves the claim.
\end{proof}
Consequently, in the situation above, $$\mr{T}^i\neq0 \quad\text{if and only if}\quad \mr{T}^{2n+1-i}\neq0.$$
On the other hand, the Hard Lefschetz isomorphism
\[
\mr{L}^{\,n-i}:\mr{H}^i(\mr{K})\xrightarrow{\sim}\mr{H}^{2n-i}(\mr{K})
\]
preserves torsion. Hence
\begin{equation*}
\mr{T}^i\xrightarrow{\ \sim\ }\mr{T}^{2n-i}
\end{equation*}
and in particular, $$\mr{T}^i\neq0\quad\text{if and only if}\quad \mr{T}^{2n-i}\neq0.$$
\begin{proposition}
\label[proposition]{torsstep}
Under the hypotheses above,
\[
\mr{T}^i\neq0\quad\Longrightarrow\quad \mr{T}^{i+1}\neq0.
\]
Similarly,
\[
\mr{T}^i\neq0\quad\Longrightarrow\quad \mr{T}^{i-1}\neq0.
\]
\end{proposition}

\begin{proof}
Suppose $\mr{T}^i\neq0$. Hard Lefschetz gives
\[
\mr{T}^{2n-i}\neq0.
\]
Applying duality to the degree $2n-i$ gives
\[
\mr{T}^{2n+1-(2n-i)}=\mr{T}^{i+1}\neq0.
\]
This proves the first implication.\\
$\hspace*{3mm}$For the second claim, using the duality, we get
\[
\mr{T}^{2n+1-i}\neq0.
\]
Applying Hard Lefschetz to the degree $i-1$, identifies torsion in degree
$2n+1-i$ with torsion in degree $i-1$. Hence $\mr{T}^{i-1}\neq0.$
\end{proof}
The first implication already suffices for the main theorem.
\2
\begin{theorem}
\label{HL-duality-flatness}
Let $R$ be a principal ideal domain and let
$\mr{K}\in \mathrm{Perf}(R)$. Assume that
\[
\mr{K}\simeq \mr{R}\operatorname{Hom}_R(\mr{K},R)[-2n]
\]
and that a degree-two operator $\mr{L}$ satisfies Hard Lefschetz
\[
\mr{L}^{\,n-i}:\mr{H}^i(\mr{K})\xrightarrow{\sim}\mr{H}^{2n-i}(\mr{K})
\]
for all $i\leq n$. Then every $\mr{H}^i(\mr{K})$ is torsion-free. Hence every
$\mr{H}^i(\mr{K})$ is a finite free $R$-module.
\end{theorem}

\begin{proof}
If $\mr{T}^i\neq0$ for some $i$,
\cref{torsstep} implies successively that
\[
\mr{T}^{i+1}\neq0,\qquad \mr{T}^{i+2}\neq0,\qquad \mr{T}^{i+3}\neq0,\quad\ldots.
\]
But $\mr{K}$ is bounded, so $\mr{H}^j(\mr{K})=0$ for $j\gg0$. This is a contradiction.
Therefore $T^i=0$ for every $i$. Since $R$ is a PID, torsion-free is equivalent to free.
\end{proof}
\2
\begin{remark}The essence of the argument is that Hard Lefschetz and duality have different symmetry centres, and their
composition is translation by one.
\end{remark}

\subsection{Application to Rees complexes}

Let $(\mr{C},\mr{F})$ be a bounded filtered complex over a field $k$, and let $\mr R_F(C)$ be its Rees complex. Assume that it is perfect over $k[\hbar]$, that it is self-dual:
\[
\mr{R}_\mr{F}(\mr{C})\simeq
\mr{R}\operatorname{Hom}_{k[\hbar]}
\left(\mr{R}_\mr{F}(\mr{C}),k[\hbar]\right)[-2n],
\]
and that it carries a degree-two operator $\mr{L}$ satisfying the Hard
Lefschetz.
Then Theorem \ref{HL-duality-flatness} shows that
$$\mr{H}\left(\mr{R}_\mr{F}(\mr{C})\right)
\text{ is free over }k[\hbar].$$
\begin{corollary}\label[corollary]{lefrees}
Let $(\mr{C,F})$ be a filtered complex whose Rees complex $\mr{R}_\mr{F}(\mr{C})\in \mr{Perf}(k[\hbar])$ satisfies Serre duality 
\[
\mr{R}_\mr{F}(\mr{C})\simeq \mr{R}\operatorname{Hom}_R\left(\mr{R}_\mr{F}(\mr{C}),R\right)[-2n]
\]
and admits a degree-two operator $\mr{L}$ with the Hard Lefschetz property
\[
\mr{L}^{\,n-i}:\mr{H}^i\left(\mr{R}_\mr{F}(\mr{C})\right)\xrightarrow{\sim}\mr{H}^{2n-i}\left(\mr{R}_\mr{F}(\mr{C})\right)
\]
for all $i\leq n$. Then the spectral sequence $$\mr{E}_1(\mr{C,F}) \Rightarrow \mr{E}_\infty(\mr{C,F})$$ degenerates on $\mr{E}_1$.
\end{corollary}

Equivalently, cohomology commutes with the classical limit:
\[
\mr{H}(\mr{R}_\mr{F}(\mr{C}))/\hbar
\xrightarrow{\ \sim\ }
\mr{H}(\mr{R}_\mr{F}(\mr{C})/\hbar)
=
\mr{H}(\operatorname{Gr}_\mr{F}(\mr{C})).
\]
\begin{example}
One way Rees complexes as in \cref{lefrees} arise is via the underlying filtered complex. Indeed, to get a perfect Rees complex, we may assume $\mr{C}$ is bounded biregular and $\operatorname{dim}_k(\mr{C}^i)<\infty$ for all $i$. Assuming there is a filtered perfect pairing $$\mr{C}\otimes\mr{C}\to k\langle w\rangle[-2n]$$ of weight $w$, applying Rees we get Serre duality (twisted by $w$) for the Rees complex. The Lefschetz condition is subtler. Even if we assume Hard Lefschetz on $\mr{C}$, generally we only get Hard Lefschetz on the central fibre of the Rees complex and Nakayama does not apply since cohomology is not a priori flat/free. In fact, the latter is precisely our goal.  
\end{example}

\section{Formality of homotopy algebras after Melani-Rubió}

\subsection{Homotopy algebras over operads and the weight filtration}

$\hspace*{3mm}$We refer the reader to \cite{loday2012algebraic} for background on operads. Let $\mathcal{P}$ be a reduced dg Koszul operad of finite type and let $\mathcal{C} \coloneqq \mathcal{P}^{\mbox{\textexclamdown}}$ be the Koszul dual cooperad. 
\2
\begin{definition}
A homotopy $\mathcal{P}$-algebra is an algebra over the Koszul dual operad $\Omega\mathcal{C}$.
\end{definition} 

$\hspace*{3mm}$A homotopy $\mathcal{P}$-algebra structure on a dg vector space $\mr{V}$ is a square-zero coderivation of degree $1$ on the cofree $\mathcal{C}$-coalgebra $\mathcal{C}(\mr{A})$ cogenerated by $\mr{A}[1]$.\\
$\hspace*{3mm}$The weight grading on $\mathcal{C}$ induces a weight decomposition 
\[\mathfrak g:=\operatorname{Coder}(\mathcal C(\mr{A}))=\prod_{w\ge 0}\mathfrak g_{(w)},\]
compatible with the graded Lie bracket:
\[[\mathfrak g_{(i)},\mathfrak g_{(j)}]\subseteq\mathfrak g_{(i+j)}.\]
$\hspace*{3mm}$In particular, if $$\mr{Q} \in \mr{Coder}^{1}(\mathcal{C}(\mr{A}))$$ is a homotopy $\mathcal{P}$-algebra structure on $\mr{A}$, then it decomposes as $$\mr{Q}=\mr{q}_1+\mr{q}_2+\cdots,$$ where each $\mr{q}_i \in \mathfrak g_{(i)}$.
We endow $\mathfrak{g}$ with a differential $$ \mr{D}_\mr{Q} = [\mr{Q},-].$$ 
$\hspace*{3mm}$The decomposition of $\mr{Q}$ induces a decomposition $$\mr{D}_\mr{Q} = \mr{D}_1+\mr{D}_2+\cdots,$$ where $\mr{D}_i = [\mr{q}_i,-]$. In particular $$\mr{D}_i : \mathfrak{g}_{(w)}^n \to \mathfrak{g}_{(w+i)}^{n+1}.$$ 
$\hspace*{3mm}$There is a complete and exhaustive weight filtration on $\mathfrak{g}$ is defined as \begin{equation}
\mr{F}^p\mathfrak g
=
\prod_{w\ge p}\mathfrak g_{(w)}.
\label{operadic-weight-filtration}
\end{equation}
Since $\mr{D}_\mr{Q}$ increases weight, we have $$\mr{D}_\mr{Q}\left(\mr{F}^p\mathfrak{g}\right) \subseteq \mr{F}^{p+1}\mathfrak{g}.$$
\begin{definition}
The operadic cohomology of the homotopy $\mathcal{P}$-algebra $\mr{A}$ is the cohomology of the deformation complex $\left(\mathfrak{g},\mr{D}_\mr{Q}\right)$.
\end{definition}
\2
\begin{example}
If $\mathcal{P}=\mr{Ass}$, we get Hochschild cohomology and for $\mathcal{P}=\mr{Lie}$, we obtain Cartan-Eilenberg cohomology.
\end{example}
\2
\begin{definition}
The operadic cohomology spectral sequence of a homotopy $\mathcal{P}$-algebra is the spectral sequence $(\mr{E}(\mr{A})_r,\mr{d}_r)$ associated to the filtered complex $(\mathfrak{g},\mr{F},\mr{D}_\mr{Q})$.
\end{definition}
 $\hspace*{3mm}$The $0^{\text{th}}$-page is just the associated graded, i.e. $$\mr{E}(\mr{A})_0^{p,q} = \mr{Hom}^{p+q}\left(\mathcal{C}(\mr{A})_{(p)},\mr{A}\right).$$ Hence the first page is $$\mr{E}(\mr{A})_1^{p,q} = \mr{Hom}^{p+q}\left(\mathcal{C}(\mr{H(A)})_{(p)},\mr{H(A)}\right)$$ and the description of $\mr{d}_1$ in terms of the connecting morphism gives $$\mr{E}(\mr{A})_2 = \mr{H}\left(\mr{Coder}(\mr{H}(\mr{A})),\mr{D}_1\right).$$ 
$\hspace*{3mm}$We are ready to state the main theorem of Melani-Robi\'o \cite{MelaniRubio}:
\2
\begin{theorem}[Melani-Rubió]\label{ml0}
Let $\mr{A}$ be a homotopy $\mathcal{P}$-algebra. The following are equivalent:
\begin{enumerate}
\item $\mr{A}$ is formal.
\item The operadic spectral sequence of $\mr{A}$ degenerates at $\mr{E}_2$.
\end{enumerate}
\end{theorem}
\2
\begin{remark}
Since the spectral sequence is homotopy invariant \cite{MelaniRubio}, we may and shall assume from now on that the homotopy $\mathcal{P}$-algebra $\mr{A}$ is minimal. 
\end{remark}
\subsection{The décalaged Rees complex}

$\hspace*{3mm}$Let's start by recalling Deligne's décalage. For a filtered complex $\mr{(C,F,d)}$, we define the décalage filtration $$\left(\mr{Dec}\mr{F}\right)^p\mr{C}^n = \{x \in \mr{F}^{p+n}\mr{C}^n \; : \; \mr{d}x\in \mr{F}^{p+n+1}\mr{C}^{n+1}\}.$$ The main property of $\mr{Dec}$ we need is the following:
\2
\begin{proposition}[Deligne]
Let $\mr{(C,F,d)}$ be a filtered complex. There is a natural isomorphism of spectral sequences $$\mr{E}_r^{p,q}(\mr{C},\mr{DecF}) \simeq \mr{E}_{r+1}^{2p+q,-p}(\mr{C},\mr{F}).$$ 
\end{proposition}
$\hspace*{3mm}$In particular, applied to the weight filtration $\mr{F}$ on $\mr{Coder}(\mathcal{C}(\mr{A}))$, it turns the theorem of Melani-Rubi\'o into an $\mr{E}_1$-degeneration statement. Let's denote the décalaged weight filtration by $\widetilde{\mr{F}}$. Because the filtration is infinite, the Rees complex associated to $$\left(\mr{Coder}(\mathcal{C}(\mr{A})),\widetilde{\mr{F}}\right)$$ is not the right object to work with in full generality. Instead, we are going to work with the pro and completed Rees complexes.
\2
\begin{remark}
 There is a graded isomorphism $$\mr{Coder}\left(\mathcal{C}(\mr{A})\right) \simeq \mr{Hom}\left(\mathcal{C}(\mr{A}),\mr{A}\right),$$ so the graded commutator induces a bracket on the right hand side, known as the Nijenhuis-Richardson bracket in the $\mathcal{P}=\mr{Lie}$  case.
\end{remark}
$\hspace*{3mm}$Recall the weight filtration $$\mr{F}^p\mathfrak{g}^{m}=\prod_{w\ge p} \mr{Hom}^m\left(\mathcal{C}(\mr{A})_{(w)},\mr{A}\right).$$ Since each component $\mr{D}_k= [q_k,-]$ of $\mr{D}$ raises the internal degree by at least $1$ as $k\ge 1$, we see that $$\widetilde{\mr{F}}\mathfrak{g}^m = \mr{F}^{p+m}\mathfrak{g}^m = \prod_{w\ge p+m} \mr{Hom}^m\left(\mathcal{C}(\mr{A})_{(w)},\mr{A}\right).$$ 
$\hspace*{3mm}$If we define the weight of an element $$f \in \mr{Hom}^m\left(\mathcal{C}(\mr{A})_{(w)},\mr{A}\right)$$ by $$\nu(f) = w-m,$$ then the décalage filtration with respect to $\nu$ is $$\widetilde{\mr{F}}^p\mathfrak{g} = \prod_{\nu\ge p} \mathfrak{g}_\nu,$$ with $$\mathfrak{g}_\nu =  \prod_{m} \mr{Hom}^m\left(\mathcal{C}(\mr{A})_{(\nu+m)},\mr{A}\right).$$ For $$f \in \mr{Hom}^m\left(\mathcal{C}(\mr{A})_{(w)},\mr{A}\right),$$ we have $$\nu(\mr{D}_kf)= w+k-(m+1)=\nu(f)+k-1,$$ so the degree of $\mr{D}_k$ with respect to $\widetilde{\mr{F}}$ is $k-1$.\\
$\hspace*{3mm}$Now for any $f \in \mathfrak{g}_w$, we have $f \in \widetilde{\mr{F}}^p\mathfrak{g}$ for all $p \le w$, hence we have an isomorphism $$\mr{R}_{\widetilde{\mr{F}},N}\mr{\left(\mathfrak{g}\right)} \simeq \mathfrak{g}\hat{\otimes} k[\hbar]/\hbar^{N+1} : \;\; f\cdot\hbar^{-p} \mapsto f\cdot\hbar^{w-p}.$$ Taking limits over $N$, we get a trivialisation of the completed Rees complex $$\widehat{\mr{R}}_{\widetilde{\mr{F}}}\left(\mathfrak{g}\right) \simeq \mathfrak{g}\brak.$$
$\hspace*{3mm}$As the weight degree of $\mr{D}_k$ is $k-1$, we have $\mr{D}_k\left(f\cdot\hbar^{-w}\right) = \left(\mr{D}_k(f)\cdot \hbar^{w+k-1}\right)\cdot \hbar^{k-1}$ and hence, under the above isomorphism, the Rees differential becomes $$\mr{D}_\mr{Q,\hbar} = \mr{D}_1+\mr{D}_2\cdot \hbar+\mr{D}_3\cdot \hbar^2+\cdots$$ Thus, the truncated/completed Rees complex becomes the constant formal deformation, but its differential now depends on $\hbar$. As $\mr{D}_\mr{Q} = [\mr{Q},-]$, we have $\mr{D}_\mr{Q,\hbar} = [\mr{Q}_\hbar,-]$, where $$\mr{Q}_\hbar = \mr{q}_1+\mr{q}_2\hbar +\mr{q}_3\hbar^2 +\cdots.$$ Since the Maurer-Cartan equation $$[\mr{Q,Q}]=0$$ is homogeneous with respect to weights, we have $$[\mr{Q_\hbar,Q_\hbar}]=0,$$ i.e. $\mr{Q}_\hbar$ defines a homotopy $\mathcal{P}$-algebra structure over the formal disc whose central fibre is $\mr{Q}_0 = \mr{q}_1$.
 
\subsection{Finite order lifting}

$\hspace*{3mm}$ Consider an element $x_0 \in \mr{H}\left(\mr{R}_{\widetilde{\mr{F}},0}(\mathfrak{g})\right).$ A lift of $x_0$ to the $N^{\text{th}}$-infinitesimal neighbourhood is represented by $x^{(N)}_\hbar = x_0+x_1\hbar+\cdots+x_N\hbar^N$ such that $$\mr{D}_\mr{Q,\hbar}\left(x^{(N)}_\hbar\right) \equiv 0 \pmod {\hbar^{N+1}}.$$ Expanding in $\hbar$, we get the equations $$\sum_{j=1}^{r}\mr{D}_{r-j+2}x_j=0, \;\; 0\le r \le N.$$ 
$\hspace*{3mm}$By \cref{infinitesimal}, these correspond precisely to the vanishing $\mr{d}_r\left(x_0\right)=0$, for $1\le r\le N$, of the spectral sequence associated to $\mr{\widetilde{F}}$. Hence, by \cref{ml0}, we deduce the following:
\2
\begin{theorem}[Melani-Rubió]\label{ml1}
Let $\mr{(A,Q)}$ be a minimal homotopy $\mathcal{P}$-algebra. The following are equivalent:
\begin{enumerate}
\item $\mr{A}$ is formal.
\item The décalaged operadic spectral sequence of $\mr{A}$ degenerates at $\mr{E}_1$.
\item Every class on the central fibre $\mr{H}\left(\mathfrak{g},[q_1,-]\right)$ lifts to $\mathbf{R}_{\widetilde{\mr{F}}}(\mathfrak{g})$.
\end{enumerate}
\end{theorem}

\subsection{Euler and Kaledin classes}

$\hspace*{3mm}$Let $\mr{(A,Q)}$ be a minimal homotopy $\mathcal{P}$-algebra. We consider the Euler derivation:
$$e_\mr{A} : \mr{A}\to \mr{A} \; : \; a \mapsto (\mr{deg}(a)+1)a.$$ 
$\hspace*{3mm}$For any $\alpha \in \mathfrak{g}^{m}_{(w)}$, we have $$\left[\alpha,e_\mr{A}\right] = (w-m)\alpha.$$ In particular, for the $\mr{q}_i$, we get $$\left[\mr{q}_i,e_\mr{A}\right] = (i-1)\mr{q}_i.$$ Hence $e_\mr{A}$ is a cocycle for $\mr{D}_{1}= [\mr{q}_1,-]$ and determines a class $$\mr{e}_\mr{A} \in \mr{H}(\mathfrak{g},[\mr{q}_1,-]) = \mr{E}(\mr{A})_2^{0,0}.$$ 
\begin{remark}
Putting the above together, we have \begin{equation}\label{var}\left[\mr{Q}_\hbar,e_\mr{A}\right] = \sum_{i \ge 1} (i-1)\mr{q}_i\hbar^{i-1} = \hbar\cdot \frac{\partial \mr{Q}_\hbar}{\partial \hbar}.\end{equation}So, the Euler derivation is a symmetry of the central fibre $\mr{q}_1$ of $\mr{Q}_\hbar$ and the obstruction that it remains a symmetry of the Rees deformation is precisely the infinitesimal variation of the deformation.
\end{remark}
$\hspace*{3mm}$Consider the right hand side of \cref{var} at $\hbar=1$: $$\dot{\mr{Q}} = \left.\frac{\partial \mr{Q}_\hbar}{\partial \hbar}\right|_{\hbar=1}.$$ Then the Maurer-Cartan equation gives $$\left[\mr{Q},\dot{\mr{Q}}\right]=0 \text{ and } \dot{\mr{Q}} \in \mr{F}^1\mathfrak{g}.$$
\begin{definition}
The operadic Kaledin class of $\mr{(A,Q)}$ is $$\mr{K}_\mr{A} = \left[\dot{\mr{Q}}\right] \in \mr{H}^{1}\left(\mr{F}^1\mathfrak{g},\mr{D}_\mr{Q}\right).$$
\end{definition}
\2
\begin{remark}
Observe that by \cref{var}, we have\[
[\mr{Q},e_\mr{A}]=\dot{\mr{Q}}.
\]
$\hspace*{3mm}$Thus $\dot{\mr{Q}}$ is exact in the full deformation complex, but its obvious
primitive $e_\mr{A}$ has weight zero and does not lie in $\mr{F}^1\mathfrak g$. The
Kaledin class measures precisely the failure of $\dot{\mr{Q}}$ to admit a primitive of strictly positive weight.
\end{remark}
\2
\begin{proposition}\label[proposition]{ek}
Let $\mr{(A,Q)}$ be a minimal homotopy $\mathcal{P}$-algebra. The Kaledin class $\mr{K}_\mr{A}$ is the obstruction class to lifting the Euler class $\mr{e}_\mr{A}$ to a Rees cocycle.
\end{proposition}
\begin{proof}
Indeed, if we have a lift $\mr{e}_{\mr{A},\hbar}=\mr{e}_\mr{A}+\hbar \cdot u(\hbar)$ of $\mr{e}_\mr{A}$, then the equation $$\mr{D}_{\mr{Q},\hbar}(\mr{e}_{\mr{A},\hbar})=0 \text{\;\; is\;\; }\hbar\cdot \frac{\partial \mr{Q}_\hbar}{\partial \hbar}+\hbar\cdot \mr{D}_{\mr{Q},\hbar}u=0$$ and since $\hbar$ acts injectively on the Rees complex, we conclude $$\frac{\partial \mr{Q}_\hbar}{\partial \hbar}=-\mr{D}_{\mr{Q},\hbar}u.$$ 
$\hspace*{3mm}$Note that $u$ has strictly positive (operadic) weight. Going backwards, we obtain a Euler lift from a positive-weight primitive of the Kaledin class.
\end{proof}
\begin{remark}
The finite truncations of the Kaledin class encode the successive obstructions to lifting the Euler
class. In particular, we obtain that the vanishing of the truncated Kaledin classes is equivalent to the survival of the Euler class through the corresponding pages of the operadic spectral sequence, cf. \cite{MelaniRubio}. Equivalently,
\begin{equation}
\mr{e}_\mr{A}\text{ lifts to }\mr{H}\left(\mr{R}_{\widetilde{\mr{F}},N}\left(\mathfrak{g}\right)\right)
\text{ if and only if } \mr{K}_\mr{A}/\hbar^{N+1}=0.
\label{eq:kaledin-euler-lifting}
\end{equation}
\end{remark}
$\hspace*{3mm}$We need a preliminary lemma before stating the main result relating the Kaledin class to Rees cohomology flatness.
\2
\begin{lemma}\label[lemma]{closedb}
Let $\mr{(A,Q)}$ be a minimal homotopy $\mathcal{P}$-algebra such that $\mr{dim}_k(\mr{A})<\infty$. Then the Rees cohomology $$\mr{H}\left(\mathfrak{g}\brak,\mr{D}_\mr{Q,\hbar}\right)$$ is $\hbar$-adically separated. 
\end{lemma}
\begin{proof}
Indeed, under the trivialisation $\widehat{\mr{R}}_{\widetilde{\mr{F}}}(\mathfrak{g}) \simeq \mathfrak{g}\brak,$ the Rees differential is $\mr{D}_{\mr{Q},\hbar} = \mr{D}_1+\mr{D}_2\cdot \hbar + \cdots$ and each $\mr{D}_k$ increases the weight by $k-1$. Hence the component of $\mr{D}_{\mr{Q},\hbar}(x)$ in any fixed weight and degree depends on only finitely many components of $x$, so $\mr{D}_{\mr{Q},\hbar}$ is continuous in the (product) weight topology.\\
$\hspace*{3mm}$As $\mr{A}$ is finite dimensional, the cohomological degree $i$ component of the Rees complex is linearly compact (in the product topology) being a product of finite dimensional vector spaces. Thus, the restriction of $\mr{D}_{\mr{Q},\hbar}$ has a closed image. Let $$[x] \in \bigcap_{N\ge 1} \hbar^{N}\mr{H}^{i}\left(\mathfrak{g}\brak,\mr{D}_{\mr{Q},\hbar}\right).$$ Then, for any $N$, we can find $y_N$ and a cocycle $z_N$ such that $$x = \mr{D}_{\mr{Q},\hbar}(y_N)+\hbar^Nz_N.$$ Taking the limit, we get $$\mr{D}_{\mr{Q},\hbar}(y_N) \longrightarrow x.$$ We have seen that $\mr{im}\left(\mr{D}_{\mr{Q},\hbar}\right)$ is closed, hence $x \in \mr{im}(\mr{D}_{\mr{Q},\hbar})$ and $[x]=0$.
\end{proof}
\2
\begin{theorem}\label{kaledinflat}
Let $\mr{(A,Q)}$ be a minimal homotopy $\mathcal{P}$-algebra such that $\mr{dim}_k(\mr{A)}<\infty$. The following are equivalent:
\begin{enumerate}
\item The Kaledin class $\mr{K}_\mr{A}$ vanishes.
\item The Rees cohomology $\mr{H}\left(\mathfrak{g}\brak,\mr{D}_{\mr{Q},\hbar}\right)$ is flat over $k\brak$.
\end{enumerate}
\end{theorem}
\begin{proof}
Assume $\mr{K}_\mr{A}=0$. Then by \cref{ek} the Euler class admits a formal lift $$\mr{e}_{\mr{A},\hbar} = \mr{e}_\mr{A}+\hbar\cdot u(\hbar).$$ Then $$\mr{E}_\hbar = \mr{ad}(\mr{e}_{\mr{A},\hbar})$$ satisfies the equation $$\left[\mr{D}_{\mr{Q},\hbar},\mr{E}_\hbar\right]=0,$$ hence we obtain an endomorphism of the Rees complex associated to the formal Euler class. On the other hand, $$\left[\mr{Q}_\hbar,e_\mr{A}\right] = \hbar\cdot \partial_\hbar \mr{Q}_\hbar$$ implies that the operator $$\mr{E}_\mr{A} = \mr{ad}(\mr{e}_\mr{A})+\hbar\cdot\partial_\hbar$$ commutes with the Rees differential too. We define $$\mr{N}_{\mr{A},\hbar} = \mr{E}_\mr{A}-\mr{E}_{\hbar} = \hbar\cdot\partial_\hbar - \mr{ad}(\hbar\cdot u(\hbar)).$$ It preserves the filtration and $$\mr{N}_{\mr{A},\hbar}(\hbar^{p}x) \equiv p\hbar^px\pmod{\hbar^{p+1}}.$$ 
$\hspace*{3mm}$Suppose the first non zero differential is $\mr{d}_r$. Then since $\mr{N}_{\mr{A},0}=p\cdot\mr{id}$ on $\mr{E}_0^{p,q}$ and the differentials $\mr{d}_i=0$ for $1\le i\le r-1$, we see $\mr{N}_{\mr{A},r}$ is still multiplication by $p$ on $\mr{E}_r^{p,q}$. Since $\mr{N}_{\mr{A},r}$ commutes with $\mr{d}_r$, we have $$\mr{N}_{\mr{A},r}(\mr{d}_rx) = \mr{d}_r(\mr{N}_{\mr{A},r}x) = p\mr{d}_rx.$$ However, $\mr{d}_rx \in \mr{E}_r^{p+r,q-r+1}$, so $$\mr{N}_{\mr{A},r}(\mr{d}_rx) = (p+r)\mr{d}_rx.$$ Hence, $$r\mr{d}_rx=0$$ for all $x \in \mr{E}_r^{p,q}$ implying $\mr{d}_r=0$ which contradicts the choice of $r$. We conclude that the spectral sequence degenerates on $\mr{E}_1$. By \cref{closedb}, the Rees cohomology is $\hbar$-adically separated and hence flat by \cref{separatede1}.\\
$\hspace*{3mm}$ The converse is trivial: by \cref{formale1}, the Euler class admits a formal lift, hence $\mr{K}_\mr{A}=0$.
\end{proof}
\begin{remark}
As in \cref{formale1}, flatness implies the Kaledin vanishing under no finite type assumptions on $\mr{A}$ or the underlying operad $\mathcal{P}$.
\end{remark}
$\hspace*{3mm}$Put together, \cref{ml1} and \cref{kaledinflat} give:
\2
\begin{theorem}\label{ml2}
Let $\mr{(A,Q)}$ be a minimal homotopy $\mathcal{P}$-algebra such that $\mr{dim}_k(\mr{A)}<\infty$. The following are equivalent:
\begin{enumerate}
\item $\mr{A}$ is formal.
\item The décalaged operadic spectral sequence of $\mr{A}$ degenerates at $\mr{E}_1$.
\item The operadic Kaledin class $\mr{K}_\mr{A}$ vanishes.
\item The operadic Rees cohomology $\mr{H}\left(\mathfrak{g}\brak,\mr{D}_\mr{Q,\hbar}\right)$ is flat over $k\brak$.
\end{enumerate}
\end{theorem}
 
 \section{Decompositions of double complexes}
 \label{khoqi}
 \subsection{Khovanov-Qi decomposition}
 
 $\hspace*{3mm}$We revisit the splitting of double complexes over a field into indecomposables following \cite{KhovanovQiTricomplexes}, see also \cite{StelzigDoubleComplexes}. 
 	\2
 \begin{definition}
 A bicomplex $\mr{C}$ is a family $\{\mr{C}^{p,q}\}$ together with two morphisms $$\partial:\mr{C}^{p,q}\to \mr{C}^{p+1,q} \text{ and } \bar\partial:\mr{C}^{p,q}\to\mr{C}^{p,q+1}$$ which satisfy the equations 
 $$\partial^2=0, \;\; \bar\partial^2=0, \;\; \partial\circ\bar\partial+\bar\partial\circ\partial=0.$$
 \end{definition}
 
 $\hspace*{3mm}$As observed in \cite{KhovanovQiTricomplexes}, bicomplexes are the same thing as bigraded left modules over the exterior algebra $\Lambda_2$ on the two generators $\partial$ and $\bar\partial$ which is bigraded by 
 $$\mr{deg}(\partial)=(1,0) \text{ and } \mr{deg}(\bar\partial)=(0,1).$$

\begin{example}[Simple bicomplex]
The most basic bicomplex is $\mr{S}^{p,q}$ given by a $1$-dimensional vector space in bidegree $(p,q)$ and $0$ elsewhere.
\end{example}
\2
\begin{example}[Square bicomplex]
The bicomplex corresponding to a (shifted) free $\Lambda_2$-module of rank $1$ is a square in the plane: $\mr{SQ}^{p,q}$ consists of $4$ $1$-dimensional vector spaces in bidegrees $(p,q), (p+1,q), (p,q+1), (p+1,q+1)$ with horizontal and vertical differentials given by the identity.
\end{example}
\2
\begin{example}[$\rightarrow$-zigzag]
The bicomplex $\mr{Z}^{p,q}_{\rightarrow,l}$ starts in bidegree $(p,q)$ and zigzags to the right and down. Here $l$ denotes the dimension of the total vector space of this bicomplex, hence it has $l-1$ non-zero arrows. We can represent it by a flat diagram as follows:
$$\{v_1\}\xrightarrow{\partial}\{v_2\}\xleftarrow{\bar\partial} \{v_3\} \xrightarrow{\partial} \cdots \xleftarrow{\bar\partial}\{v_l\}.$$
\end{example}
\2
\begin{example}[$\uparrow$-zigzag]
The bicomplex $\mr{Z}^{p,q}_{\uparrow,l}$ starts in bidegree $(p,q)$ and zigzags down and to the right. Here $l$ denotes the dimension of the total vector space of this bicomplex, hence it has $l-1$ non-zero arrows. This time, the flat diagram is as follows:
$$\{v_1\} \xleftarrow{\bar\partial}\{v_2\}\xrightarrow{\partial} \{v_3\} \xleftarrow{\bar\partial} \cdots \xrightarrow{\partial}\{v_l\}.$$
\end{example}

$\hspace*{3mm}$We are ready to state the decomposition theorem for (bounded) bicomplexes.
\2
\begin{theorem}[\cite{KhovanovQiTricomplexes}]
Any bounded bicomplex $\mr{C}$ decomposes into a direct sum of the indecomposable bicomplexes $\mr{S}^{p,q}, \mr{SQ}^{p,q}, \mr{Z}^{p,q}_{\rightarrow,l} \text{ and } \mr{Z}^{p,q}_{\uparrow,l}$.
\end{theorem}

\subsection{Rees description}
$\hspace*{3mm}$As the Rees functor is additive, we can apply it to the above structure theorem. Consider a bicomplex $\mr{C} = \{\mr{C}^{p,q}\}$ and its total complex $\mr{Tot}(\mr{C})$. We can filter $\mr{Tot}(\mr{C})$ by columns, i.e. $$\mr{F}^p\mr{Tot}^n(\mr{C}) = \bigoplus_{i \ge p}\mr{C}^{i,n-i}.$$
$\hspace*{3mm}$The Rees complex $\mr{R}(\mr{Tot(C)}$ is a differential graded module over $k[\hbar]$. If we treat it as a graded $k[\hbar]$-submodule of $k[\hbar,\hbar^{-1}]$, then the Rees differential is simply $$\mr{d_R}=\bar\partial+\partial.$$ 
$\hspace*{3mm}$On the other hand, to any $v \in \mr{C}^{p,q}$, we have an associated Rees basis element $v\hbar^{-p}$. The vertical differential preserves the column filtration, while the horizontal differential shifts the filtration by $+1$. Hence the equation for the Rees differential in the Rees basis is \begin{equation}\label{rd}\mr{d}_\mr{R} = \bar\partial+\hbar\partial.\end{equation} 
$\hspace*{3mm}$With this preparation, we are now going to apply the Rees functor to the indecomposable bicomplexes described in the examples above.
\2
\begin{example}[Simples]
There is nothing going on here, to a simple in the plane, sitting in bidegree $(p,q)$, we associate a rank $1$ free module $k[\hbar]$ with no differential. The spectral sequence is trivial as the cohomology is free.
\end{example}
\2
\begin{example}[Squares]
In this case, since the total degrees are $p+q, p+q+1 \text{ and } p+q+2$, according to our equation for the Rees differential, we get the complex 
$$k[\hbar] \xrightarrow{\begin{pmatrix}
1\\ \hbar
\end{pmatrix}}
k[\hbar]\oplus k[\hbar] \xrightarrow{\begin{pmatrix}
-\hbar & 1
\end{pmatrix}}
k[\hbar].$$
$\hspace*{3mm}$Further, we can change the basis in the middle to $$\left\{\begin{pmatrix}
1\\ \hbar
\end{pmatrix}, \begin{pmatrix}
0\\ 1
\end{pmatrix}\right\},$$ hence the complex splits as $$(k[\hbar] \xrightarrow{\ 1\ } k[\hbar]) \;\oplus\; (k[\hbar] \xrightarrow{\ 1\ } k[\hbar]).$$
$\hspace*{3mm}$ This, as well as its reduction$\mod\hbar$, has trivial cohomology, hence already the first page of the spectral sequence vanishes and so these indecomposables do not contribute to the spectral sequence anything at all.
\end{example}
\2
\begin{example}[Horizontal zigzags]
The horizontal zigzag, starting in bidegree $(p,q)$, has sources for the differential given by $\{v_1,v_3,\cdots,\}$ and targets $\{v_2,v_4,\cdots,\}$. The sources all have total degree $n=i+j$, while the targets get $p+q+1=n+1$. To each $v_i$, we get a homogeneous Rees element $v_i\hbar^{-p_i}$. Using the equation for the Rees differential \c	ref{rd}, we distinguish two cases.
\begin{enumerate}
\item If $l=2r$, we get a Jordan block $$\mr{d}_{\mr{R}} = \begin{pmatrix}
\hbar & 1 & 0 & \dots & 0 \\
0 & \hbar & 1 & \dots & 0 \\
0 & 0 & \hbar & \dots & 0 \\
\vdots & \vdots & \vdots & \ddots & \vdots \\
0 & 0 & 0 & \dots & \hbar
\end{pmatrix}.$$
Upon changing basis, the matrix is $\mr{diag}(1,1,\cdots,1,\hbar^r)$. Hence the Rees complex splits as $$\left(k[\hbar]\xrightarrow{\hbar^r}k[\hbar]\right)\oplus\left(k[\hbar]\xrightarrow{1}k[\hbar]\right)^{\oplus (r-1)}.$$ The cohomology is a cyclic module $k[\hbar]/\hbar^r$ concentrated in degree $n+1$, i.e. it's purely torsion and this corresponds precisely to the non-triviality of the spectral sequence which has a single non-zero differential $$\mr{E}_r^{p,q}\xrightarrow{1} \mr{E}_r^{p+r,q-r+1}.$$
\item If $l=2r+1$ is odd, the 	$r\times(r+1)$ matrix representing the Rees differential takes the form $$\mr{d}_{\mr{R}}= \begin{pmatrix}
\hbar & 1 & 0 & \dots & 0 & 0 \\
0 & \hbar & 1 & \dots & 0 & 0 \\
\vdots & \vdots & \vdots & \ddots & \vdots & \vdots \\
0 & 0 & 0 & \dots & \hbar & 1
\end{pmatrix}.$$
$\hspace*{3mm}$As above, the Rees complex splits as $$\left(k[\hbar]\xrightarrow{0} 0\right)\oplus\left(k[\hbar]\xrightarrow{1} k[\hbar]\right)^{\oplus r}.$$
$\hspace*{3mm}$The cohomology is free of rank $1$ sitting in total degree $i+j=n$. This corresponds to the fact that the spectral sequence has a single non-zero entry $\mr{E}_1^{p,q}=k$ which survives all pages. $\hspace*{3mm}$Hence each of these indecomposables contributes exactly one $1$-dimensional permanent cycle (a cycle surviving to $\mr{E}_\infty$) to the spectral sequence. 
\end{enumerate}
\end{example}
\2
\begin{example}[Vertical zigzags]
This is very similar to the previous example, in fact, as the zigzags begin with a vertical arrow, we now have sources at $\{v_2,v_4,\cdots,\}$ and targets $\{v_1,v_3,\cdots,\}$. 
\begin{enumerate}
\item For even $l=2r$, we get $$\mr{d}_{\mr{R}} = \begin{pmatrix}
1 & \hbar & 0 & \dots & 0 \\
0 & 1 & \hbar & \dots & 0 \\
0 & 0 & 1 & \dots & 0 \\
\vdots & \vdots & \vdots & \ddots & \vdots \\
0 & 0 & 0 & \dots & 1
\end{pmatrix}.$$
The Rees complex splits as $$\left(k[\hbar]\xrightarrow{1} k[\hbar]\right)^{\oplus r}.$$ The cohomology vanishes. In fact, already the first page of the spectral sequence $\mr{E}_1=0$ which precisely corresponds to the fact that Rees cohomology is trivial instead of torsion.
\item The odd case $l=2r+1$ has $(r+1)\times r$ matrix $$\mr{d}_{\mr{R}}= \begin{pmatrix}
\hbar & 0 & \dots & 0 \\
1 & \hbar & \dots & 0 \\
0 & 1 & \dots & 0 \\
\vdots & \vdots & \ddots & \vdots \\
0 & 0 & \dots & \hbar \\
0 & 0 & \dots & 1
\end{pmatrix}.$$
The Smith normal form gives a splitting $$\left(0\xrightarrow{0} k[\hbar]\right)\oplus\left(k[\hbar]\xrightarrow{1} k[\hbar]\right)^{\oplus r}.$$ The cohomology is again free of rank $1$, concentrated in total degree $p+q=n+1$. There is a single non-zero entry on the first page of the spectral sequence $\mr{E}_1^{p+r,q-r}=k$ which is a permanent cycle surviving all the way to $\mr{E}_\infty$ as in the odd length horizontal zigzags.
\end{enumerate}
\end{example}

$\hspace*{3mm}$ We can summarise the discussion in this section as follows: any (bounded) bicomplex splits into indecomposable complexes and there are $4$ types of these - simples, squares, horizontal and vertical zigzags. According to the analysis above, their relations to the Rees complex and the spectral sequence of the double complex is:
\begin{enumerate}
\item Each simple contributes a rank $1$ free summand to the cohomology of the Rees complex and hence a single permanent cycle that survives all pages of the spectral sequence.
\item Squares are invisible in cohomology, hence they do not feature in the spectral sequence as well.
\item Horizontal zigzags have two subcases. Odd length zigzags adds a rank $1$ free summand to the cohomology of the Rees complex. Hence, in terms of the spectral sequence, we get a single cycle on the first page that survives all the way to $\mr{E}_\infty$ and remains as a permanent cycle, contributing one dimension.\\ Even length zigzags are the source of torsion in the cohomology of the Rees complex. Each of length $2r$ contributes a cyclic module $k[\hbar]/\hbar^r$ to the cohomology and thus the spectral sequence is non-trivial with a single non-zero differential $\mr{d}_r:\mr{E}_r\to \mr{E}_r$ which is the identity.
\item Vertical zigzags also get two subcases. The even length vertical zigzags vanish in cohomology, and so do not alter the spectral sequence too. \\ Odd length vertical zigzags behave like the odd length horizontal zigzags adding rank $1$ free summand to Rees cohomology and hence a single permanent cycle to the spectral sequence.
\end{enumerate}

\section{The two-parameter Rees complex of a double complex}

\label{sec:two-parameter-rees-double-complex}

$\hspace*{3mm}$We now refine the one-parameter Rees description of Section~\ref{khoqi} by retaining both filtrations of a
double complex simultaneously. The resulting two-parameter object is much more rigid than either of its one-parameter restrictions. It detects not only the even zigzags responsible for the two spectral sequences, but also the odd zigzags and the square summands which are invisible to one or both one-parameter Rees cohomologies.
\subsection{Two-parameter Rees construction}
$\hspace*{3mm}$ We begin with a definition that works for any bifiltered module.
\2
\begin{definition}
The two-parameter Rees module of a bifiltered module $(\mr{M,F,G})$ is the graded $k[s,t]$-module $$\mr{R}_{F,G}(\mr{M}) = \bigoplus_{p,q} \left(\mr{F}^p\mr{M}\cap\mr{G}^q\mr{M}\right)s^{-p}t^{-q}.$$
\end{definition}
 $\hspace*{3mm}$Let $(\mr{C}^{p,q},\partial,\bar\partial)$
be a bounded double complex over a field $k$ and let
\[
R:=k[s,t],
\qquad
\mathfrak m:=(s,t).
\]
$\hspace*{3mm}$With the natural $\mathbf Z^2$-grading on $s$ and $t$, the two-parameter Rees construction for the column and row
filtrations of $\mr{Tot}(\mr{C})$ is quasi-isomorphic to
\begin{equation*}
\mr{R}_{s,t}(\mr{C})
:=
\left(
\operatorname{Tot}(\mr{C})\otimes_k R,
\mr{d}_{s,t}
\right),
\qquad
\mr{d}_{s,t}=s\partial+t\bar\partial.
\end{equation*}
\begin{proposition}
We have the standard relations:
\begin{align*}
\mr{R}_{s,t}(\mr{C})/(t-1)
&=
\left(
\operatorname{Tot}(\mr{C})[s],
\bar\partial+s\partial
\right),
\\
\mr{R}_{s,t}(\mr{C})/(s-1)
&=
\left(
\operatorname{Tot}(\mr{C})[t],
\partial+t\bar\partial
\right),
\\
\mr{R}_{s,t}(\mr{C})/(s-1,t-1)
&=
\left(
\operatorname{Tot}(\mr{C}),
\partial+\bar\partial
\right),
\\
\mr{R}_{s,t}(\mr{C})/(s,t)
&=
\left(
\operatorname{Tot}(\mr{C}),
0
\right).
\end{align*}
\end{proposition}

\subsection{Indecomposables and the two-parameter Rees}
$\hspace*{3mm}$Before computing the two-parameter Rees of the Khovanov-Qi indecomposables, it would be useful to set up some notation. We shall repeatedly need the matrices $$J_n^{s}(s,t) \text{ and } J_n^{t}(s,t) \in M_n(R)$$ defined by
\begin{equation*}
J_n^{s}(s,t)
:=
\begin{pmatrix}
s&t&0&\cdots&0\\
0&s&t&\ddots&\vdots\\
\vdots&\ddots&\ddots&\ddots&0\\
0&\cdots&0&s&t\\
0&\cdots&\cdots&0&s
\end{pmatrix},
\end{equation*}
and
\begin{equation*}
J_n^{t}(s,t)
:=
\begin{pmatrix}
t&s&0&\cdots&0\\
0&t&s&\ddots&\vdots\\
\vdots&\ddots&\ddots&\ddots&0\\
0&\cdots&0&t&s\\
0&\cdots&\cdots&0&t
\end{pmatrix}.
\end{equation*}
$\hspace*{3mm}$Note that particular instances of these already appeared in the previous section as differentials of the Rees complexes. We write
\begin{equation}
Q_n^s:=\operatorname{coker}J_n^s,
\qquad
Q_n^t:=\operatorname{coker}J_n^t.
\end{equation}
 for the cokernels of the above matrices.\\
$\hspace*{3mm}$For the odd zigzags we shall need
\begin{equation*}
\mr{A}_r(s,t)
:=
\begin{pmatrix}
s&t&0&\cdots&0\\
0&s&t&\ddots&\vdots\\
\vdots&\ddots&\ddots&\ddots&0\\
0&\cdots&0&s&t
\end{pmatrix}
:
R^{r+1}\longrightarrow R^r.
\end{equation*}
and the vector
\begin{equation*}
v_r
=
\begin{pmatrix}
t^r\\
-st^{r-1}\\
s^2t^{r-2}\\
\vdots\\
(-s)^r
\end{pmatrix}
\end{equation*}
which is in the kernel of $\mr{A}_r(s,t)$:
\[
\mr{A}_rv_r=0.
\]

$\hspace*{3mm}$Finally, define
\[
\mathfrak m^r=(s,t)^r
\]
and let
\begin{equation*}
\omega_r
:=
\operatorname{Ext}^2_R(R/\mathfrak m^r,R)
\end{equation*}
be the canonical module of the Artinian ring $R/\mathfrak m^r$.
\2
\begin{remark}
We shall suppress all grading and cohomological shifts until they are
needed. Signs in the matrices depend on the convention
\[
\partial\bar\partial+\bar\partial\partial=0;
\]
they may be removed by changes of basis and do not affect the module
calculations.
\2
\end{remark}
\begin{example}[Simples]

For a simple $\mr{S}^{p,q}$, both differentials vanish. Therefore
\begin{equation}
\mr{R}_{s,t}(\mr{S}^{p,q})
\simeq
R
\end{equation}
up to the corresponding internal and total-degree shift. So, simples give free rank-one modules supported on
\[
\mathbf A^2_{s,t}.
\]
\end{example}
\2
\begin{example}[Squares]

Let $\mr{SQ}$ be a square. After choosing bases, its two-parameter totalisation is the Koszul complex
\begin{equation*}
0
\longrightarrow
R
\xrightarrow{\binom{s}{t}}
R^2
\xrightarrow{\begin{pmatrix}-t&s\end{pmatrix}}
R
\longrightarrow
0.
\end{equation*}
Since $s,t$ form a regular sequence, its only cohomology is
\begin{equation*}
\mr{H}(\mr{R}_{s,t}(\mr{SQ}))
\simeq
R/(s,t)
=
k
\end{equation*}
in the top cohomological degree.\\
$\hspace*{3mm}$Hence a square, which is invisible to either one-parameter Rees/spectral sequence as we saw in the previous section, becomes a codimension-two skyscraper at the origin in the two-parameter family. Notice that
\[
\mr{K}(s,t)\otimes_RR/(t-1)
\]
and
\[
\mr{K}(s,t)\otimes_RR/(s-1)
\]
are contractible, because one of the two Koszul generators becomes a unit.
This recovers the fact that a square contributes nothing to either
one-parameter Rees cohomology.
\end{example}
\2
\begin{example}[Even zigzags]
We first consider a zigzag with an odd number of arrows,
\[
2n-1,
\]
and hence $2n$ vertices.\\
$\hspace*{3mm}$We've seen such a zigzag is responsible for a nonzero spectral sequence differential. For a horizontal zigzag, i.e. the orientation beginning with a horizontal arrow, the two adjacent
total degrees each have rank $n$, and the two-parameter differential is
\begin{equation*}
R^n
\xrightarrow{J_n^s(s,t)}
R^n.
\end{equation*}
Since
$
\det J_n^s=s^n,
$
the map is injective and
\begin{equation*}
\mr{H}(\mr{R}_{s,t}(Z_{\rightarrow,2n-1}))
\simeq
Q_n^s.
\end{equation*}
Similarly, the opposite orientation, corresponding to a vertical zigzag, gives
\begin{equation*}
\mr{H}(\mr{R}_{s,t}(Z_{\uparrow,2n-1}))
\simeq
Q_n^t.
\end{equation*}
\end{example}
\begin{remark}
The zeroth Fitting ideals are
\begin{equation*}
\operatorname{Fitt}_0(Q_n^s)=(s^n),
\qquad
\operatorname{Fitt}_0(Q_n^t)=(t^n).
\end{equation*}
Hence $Q_n^s$ is supported on the $n$-fold thickening of the
$t$-axis given by $s^n$,
whereas $Q_n^t$ is supported on the thickening of the
$s$-axis given by $t^n$. After localising away from the origin along the corresponding axis, we get
\begin{equation*}
(Q_n^s)_t
\simeq
R_t/(s^n),
\qquad
(Q_n^t)_s
\simeq
R_s/(t^n).
\end{equation*}
\end{remark}
\2
\begin{remark}
The one-parameter restrictions recover the familiar torsion blocks:
\begin{align*}
Q_n^s\otimes_RR/(t-1)
&\simeq
k[s]/(s^n),
\\
Q_n^s\otimes_RR/(s-1)
&=0,
\\
Q_n^t\otimes_RR/(s-1)
&\simeq
k[t]/(t^n),
\\
Q_n^t\otimes_RR/(t-1)
&=0.
\end{align*}
\end{remark}

\begin{example}[Odd zigzags]
We now consider a zigzag with $2r$ arrows and $2r+1$ vertices. Such a
zigzag contributes one total cohomology class and is invisible to the higher differentials of both spectral sequences. Nevertheless, the two-parameter Rees complex detects its length.
\begin{itemize}
\item For a vertical zigzag, its total complex has $r$ generators in the lower degree and $r+1$ in the upper degree, the two-parameter complex is
\begin{equation*}
0
\longrightarrow
R^r
\xrightarrow{A_r(s,t)^T}
R^{r+1}
\longrightarrow
0.
\end{equation*}

We have a resolution of the ideal $\mathfrak m^r$ given by
\begin{equation}
0
\longrightarrow
R^r
\xrightarrow{A_r^T}
R^{r+1}
\xrightarrow{v_r^T}
\mathfrak m^r
\longrightarrow
0,
\label{hb}
\end{equation}
where the components of $v_r$ are the degree-$r$ monomials
\[
t^r,-st^{r-1},\ldots,(-s)^r.
\]
Therefore
\begin{equation*}
\mr{H}\left(\mr{R}_{s,t}\left(\mr{Z}_{\uparrow,2r}\right)\right)
\simeq
\mathfrak m^r.
\end{equation*}
The ideal $\mathfrak m^r$ is a torsion-free module of rank one. Away from the origin, we have
\[
\mathfrak m^r|_{\mathbf A^2\setminus\{0\}}
\simeq
\mathscr O_{\mathbf A^2\setminus\{0\}},
\]
so the odd zigzag looks like a simple on the punctured plane.
It differs from a simple at the origin:
\begin{equation*}
(\mathfrak m^r)^{\vee\vee}/\mathfrak m^r
\simeq
R/\mathfrak m^r.
\end{equation*}
The length of the latter is
\begin{equation*}
\operatorname{length}(R/\mathfrak m^r)
=
\frac{r(r+1)}2.
\end{equation*}

\item For the opposite orientation, the total complex is
\begin{equation}
0
\longrightarrow
R^{r+1}
\xrightarrow{A_r(s,t)}
R^r
\longrightarrow
0.
\end{equation}

Dualising \cref{hb} gives
\begin{equation*}
0
\longrightarrow
\operatorname{Hom}_R(\mathfrak m^r,R)
\longrightarrow
R^{r+1}
\xrightarrow{A_r}
R^r
\longrightarrow
\operatorname{Ext}^1_R(\mathfrak m^r,R)
\longrightarrow0.
\label{dual-hb}
\end{equation*}
Since $R$ is normal and $\mathfrak m^r$ has reflexive hull $R$,
\[
\operatorname{Hom}_R(\mathfrak m^r,R)\simeq R,
\]
while from
\[
0\to\mathfrak m^r\to R\to R/\mathfrak m^r\to0
\]
one obtains
\[
\operatorname{Ext}^1_R(\mathfrak m^r,R)
\simeq
\operatorname{Ext}^2_R(R/\mathfrak m^r,R)
=
\omega_r.
\]
Hence the Rees cohomology is
\begin{equation*}
\mr{H}^{*}(\mr{R}_{s,t}(\mr{Z}_{\rightarrow,2r}))
\simeq
R,
\qquad
\mr{H}^{*+1}(\mr{R}_{s,t}(\mr{Z}_{\rightarrow,2r}))
\simeq
\omega_r
\end{equation*}
for the two non-zero total degrees. 
\end{itemize}
\end{example}
\2
\begin{remark}
We can rephrase the result more conceptually as
\begin{equation*}
\mr{R}_{s,t}(\mr{Z}_{\rightarrow,2r})
\simeq
\mr{RHom}_R(\mathfrak m^r,R)
\end{equation*}
up to the appropriate shift. Thus the two orientations of an odd zigzag are exchanged by derived duality.
\end{remark}
\2
\begin{remark}
The distinction between a simple and an odd zigzag disappears on either affine line $t=1$ or $s=1$ since
\begin{equation*}
\mathfrak m^r\otimes_RR/(t-1)
\simeq
k[s] \text{\quad and \quad}
\mathfrak m^r\otimes_RR/(s-1)
\simeq
k[t].
\end{equation*}
$\hspace*{3mm}$For the dual orientation, the finite-length module $\omega_r$ is supported at the origin and therefore disappears under either specialization
$t=1$ or $s=1$, so only the free $R$-summand remains.
$\hspace*{3mm}$Hence both one-parameter spectral sequences see only a permanent cycle, while the two-parameter construction remembers the length of the odd zigzag through the
singularity at
\[
(s,t)=(0,0).
\]
This is one reason why one needs more than degeneration of both
spectral sequences in order to obtain the
$\partial\bar\partial$-lemma.\\
$\hspace*{3mm}$Also, note that, keeping the natural $\mathbf Z^2$-grading, the variables $s$ and
$t$ remember which arrows were horizontal and which were vertical.
\end{remark}

\subsection{Bott-Chern and Aeppli cohomology}

$\hspace*{3mm}$The two-parameter calculation gives a direct algebraic interpretation of the Bott-Chern and Aeppli contributions of a zigzag. For a zigzag, we have $$\partial\bar\partial=0$$ and so
\2
\begin{itemize}
\item every vertex with no outgoing arrow contributes one
Bott-Chern class;
\item every vertex with no incoming arrow contributes one
Aeppli class.
\end{itemize}
These are respectively the sinks and sources of the oriented zigzag, which is the square--zigzag description of Bott-Chern and Aeppli
cohomology of \cite{StelzigDoubleComplexes}.\\
$\hspace*{3mm}$For an even vertex differential zigzag there are equally many sources and
sinks:
\begin{equation*}
\dim \mr{H}_\mr{BC}
=
\dim \mr{H_A}
=
n,
\qquad
\text{for a $2n$-vertex zigzag}.
\end{equation*}
$\hspace*{3mm}$For an odd vertex zigzag of dimension $2r+1$ represented by
\[
R^r\to R^{r+1},
\]
there are $r$ sources and $r+1$ sinks, so
\begin{equation*}
\dim \mr{H_A}=r,
\qquad
\dim \mr{H_{BC}}=r+1.
\end{equation*}
The opposite orientation is represented by
\[
R^{r+1}\to R^r,
\]
and we get the reversed numbers:
\begin{equation*}
\dim \mr{H_A}=r+1,
\qquad
\dim \mr{H_{BC}}=r.
\end{equation*}
$\hspace*{3mm}$A simple contributes one class to both and a square contributes none, since its unique common kernel is the image of $\partial\bar\partial$. Thus the Bott-Chern and Aeppli cohomology dimensions are local invariants of the two-parameter Rees complex at the origin.

\subsection{The \texorpdfstring{$\partial\bar\partial$}{ddbar}-lemma via the Rees complex}

$\hspace*{3mm}$A theorem of Stelzig in \cite{StelzigDoubleComplexes} says that a double complex satisfies the
$\partial\bar\partial$-lemma precisely when the decomposition of the double complex
contains only simples and squares:
\2
\begin{theorem}[Stelzig]\label[proposition]{stelzig}
A bounded double complex satisfies the
$\partial\bar\partial$-lemma if and only if its decomposition into indecomposables contains only simples and squares.
\end{theorem}
\2
\begin{remark}
 The Rees reformulation of \cref{stelzig} goes as follows. Define
\[
\mathcal D_{\partial\bar\partial} \subseteq \operatorname{Perf}^{\mathbf Z^2}(R)
\]
to be the smallest additive idempotent-complete subcategory generated by
all bigrading and cohomological shifts of
\[
R
\qquad\text{and}\qquad
k=R/\mathfrak{m}.
\]

Then \cref{ddctwo} may be restated as
$$
\mr{C}\text{ satisfies the }\partial\bar\partial\text{-lemma}
\text{ if and only if }
\mr{R}_{s,t}(\mr{C})\in\mathcal D_{\partial\bar\partial}.$$
\end{remark}
$\hspace*{3mm}$The decomposition of a double complex into indecomposables is not canonical and we would like to formulate an intrinsic version of this result which doesn't require choosing such a decomposition.\\
Squares satisfy the $\partial\bar{\partial}$-lemma, but their two-parameter Rees complex is the Koszul complex $\mr{K}_R(s,t)$ whose cohomology is not flat. Hence requiring flatness over $\mathbf{A}^2$ is too strong. In fact, we have the following characterisation.
\2
\begin{proposition}
Let $(\mr{C},\partial,\bar\partial)$ be a bounded double complex of finite dimensional vector spaces. The following are equivalent:
\begin{enumerate}
\item Every $\mr{H}^m(\mr{R}_{s,t}(\mr{C}))$ is flat over $R$.
\item Every $\mr{H}^m(\mr{R}_{s,t}(\mr{C}))$ is free over $R$.
\item The decomposition of $\mr{C}$ contains only simples.
\item $\mr{R}_{s,t}(\mr{C})$ is a direct sum of bigraded shifts of $R$.
\end{enumerate}
\end{proposition}
\begin{proof}
Since $R=k[s,t]$ is Noetherian and the cohomology modules are finitely
generated, flatness is equivalent to projectivity, and projective modules
over $R$ are free.\\
$\hspace*{3mm}$The elementary calculations of
Section~\ref{sec:two-parameter-rees-double-complex} show that every
non-simple indecomposable contributes a non-flat module:
\begin{itemize}
\item a square gives $R/\mathfrak m$;
\item an even zigzag contributes a torsion module $Q_n^s$ or
$Q_n^t$;
\item a vertical odd zigzag contributes
$\mathfrak m^r$, which is torsion-free but not free for $r\ge1$;
\item the horizontally orientated odd zigzags give a finite-length module in addition
to a free summand.
\end{itemize}
$\hspace*{3mm}$Thus flat cohomology occurs precisely for sums of simples.
\end{proof}

$\hspace*{3mm}$Denote
$\mr{U}:=\mathbf A^2\setminus\{0\}$ and let
$j:\mr{U}\hookrightarrow \mathbf A^2$
be the inclusion. Working over the punctured plane instead, we eliminate only the even zigzags:
\2
\begin{proposition}[Flatness over the punctured plane]
Let $(\mr{C},\partial,\bar\partial)$ be a bounded double complex. The following are equivalent:
\begin{enumerate}
\item $j^*\mr{H}\left(\mr{R}_{s,t}(\mr{C})\right)$ is locally free on $\mr{U}$.
\item The decomposition of $\mr{C}$ contains no even
zigzags.
\item Both spectral sequences associated with the two filtrations of the
double complex degenerate at $E_1$.
\end{enumerate}
\end{proposition}
\begin{proof}
Simples give free modules on $\mathbf{A}^2$, while squares are supported only at the
origin and therefore vanish on $\mr{U}$.\\
$\hspace*{3mm}$An even zigzag of one orientation contributes $Q_n^s$, which is supported
on the thickened axis
\[
V(s^n),
\]
and is nonzero at every generic point of the $t$-axis. Thus it is not
locally free as an $R$-module on $\mr{U}$. The opposite orientation gives the
same thing along $V(t^n)$ instead.\\
$\hspace*{3mm}$On the other hand, an odd zigzag contributes either $\mathfrak m^r$ or its
derived dual. Since
\[
\mathfrak m^r|_\mr{U}\simeq\mathscr O_\mr{U}
\]
and the finite-length cohomology of the dual orientation is supported only at
the origin, every odd zigzag has locally free cohomology on $\mr{U}$.\\
$\hspace*{3mm}$Hence flatness over $\mr{U}$ is equivalent precisely to the absence of
even zigzags. By the square-zigzag description of the two spectral
sequences, this is equivalent to $E_1$-degeneration in both directions, see e.g.
\cite{KhovanovQiTricomplexes,StelzigDoubleComplexes}.
\end{proof}
$\hspace*{3mm}$So, imposing flatness over $\mr{U}$ is equivalent to having no even zigzags. Next, we need a condition that, paired with flatness over $\mr{U}$, is equivalent to having no odd zigzags which will leave us with simples and squares that is precisely the $\partial\bar\partial$-lemma according to \cref{stelzig}. Looking at the Rees complex of an odd zigzag, we have three independent ways it can appear at the origin. \\
$\hspace*{3mm}$Let
\[
\mr{M}^m:=\mr{H}^m\left(\mr{R}_{s,t}(\mr{C})\right)
\]
and
\begin{equation*}
\mr{T}^m
:=
H^0_{\mathfrak m}(\mr{M}^m)
=
\{x\in \mr{M}^m:\mathfrak m^N x=0\text{ for some }N\}
\end{equation*}
be the maximal finite-length submodule. Finally, put
\[
\overline{\mr{M}^m}:=\mr{M}^m/\mr{T}^m.
\]
$\hspace*{3mm}$Recall that for an odd vertical zigzag with $2r+1$ vertices, we have only one non-zero cohomology which is $$\mr{M}\simeq \mathfrak{m}^r.$$ This is a torsion-free rank $1$ module over $R$ but it is not reflexive and $$\left(\mathfrak{m}\right)^{\vee\vee}\simeq R \text{ and } \left(\mathfrak{m}\right)^{\vee\vee}/\mathfrak{m} \simeq R/\mathfrak{m}.
$$ 
$\hspace*{3mm}$For the opposite orientation, the additional finite-length cohomology is
\[
\omega_r
=
\operatorname{Ext}^2_R(R/\mathfrak m^r,R).
\]
This is the canonical module of the Artinian ring $R/\mathfrak m^r$.
Its annihilator is
\[
\operatorname{Ann}_R(\omega_r)=\mathfrak m^r.
\]
and so
$$\mathfrak m\omega_r=0 \iff r=1.$$
$\hspace*{3mm}$Finally, for a horizontal zigzag, when $r=1$, we have three vertices and the complex has cohomology $$\mr{H}^0 \simeq R \text{ and } \mr{H}^1 \simeq k,$$ up to shifts, and hence sits in a triangle $$R\longrightarrow \mr{R}_{s,t}(\mr{C})\longrightarrow k[-1]\xrightarrow{\;\;\epsilon\;\;} R[1].$$ 
$\hspace*{3mm}$The Postnikov class $$\epsilon \in \mr{Ext}^2(k,R)$$ is what distinguishes the three-vertex zigzag from a sum of a simple and a square. 
\2
\begin{definition}
We say that $\mr{R}_{s,t}(\mr{C})$ is \emph{cohomologically reflexive} if, for every
$m$,
$$
\overline{\mr{M}}^m
\longrightarrow
(\overline{\mr{M}}^m)^{\vee\vee}$$
is an isomorphism.
\end{definition}
$\hspace*{3mm}$If punctured-plane flatness already holds, being cohomologically reflexive  removes exactly the odd zigzags whose cohomology is $\mathfrak m^r$.
\2
\begin{remark}
Suppose $\mr{M}$ is a torsion-free $R$-module and $$\widetilde{\mr{M}}|_{\mr{U}}$$ is flat, i.e. the conditions the cohomology modules of the Rees complex we require satisfy. We have a natural map $$\widetilde{\mr{M}}\longrightarrow j_*\left(\widetilde{\mr{M}|_\mr{U}}\right).$$ 
Since $\mr{R}$ is regular, $$
j_*\left(\widetilde{\mr{M}|_\mr{U}}\right)\simeq \widetilde{\mr{M}^{\vee\vee}}.
$$ Hence $\mr{M}$ is reflexive iff it satisfies the Hartogs' condition $$\widetilde{\mr{M}}\simeq j_*\left(\widetilde{\mr{M}|_\mr{U}}\right).$$
\end{remark}
\2
\begin{definition}
We say that $\mr{R}_{s,t}(\mr{C})$ has \emph{reduced origin torsion} if
$$\mathfrak m \mr{T}^m=0$$
for every $m$.
\end{definition}
$\hspace*{3mm}$Requiring reduced torsion on top of the previous two conditions removes all odd horizontal zigzags for $r\ge 2$. We are thus left with simples, squares and three-vertex horizontal zigzags. To remove the latter, we require that the Postnikov $\mr{Ext}^2$-classes vanish.
\2
\begin{remark}
Having imposed the first three conditions -- flatness over $\mr{U}$, cohomological reflexivity and reduced torsion at the origin, the vanishing of the Postnikov classes is equivalent to the formality of the complex $\mr{R}_{s,t}(\mr{C})$. 
\end{remark}
$\hspace*{3mm}$We can summarise the discussion as follows:
\2
\begin{theorem}[Two-parameter Rees criterion for the
$\partial\bar\partial$-lemma]
\label{ddctwo}
Let $\mr{C}$ be a bounded double complex of finite-dimensional vector spaces and
let
\[
\mr{R}_{s,t}(\mr{C})
\]
be its two-parameter Rees complex over $R=k[s,t]$.
The following are equivalent.

\begin{enumerate}
\item $\mr{C}$ satisfies the $\partial\bar\partial$-lemma.

\item In the decomposition of $\mr{C}$, only simples and squares
occur.

\item The following intrinsic conditions hold:
\begin{enumerate}
\item[(a)] $\mr{H}\left(\mr{R}_{s,t}(\mr{C})\right)$ is locally free on
\[
\mathbf A^2\setminus\{0\};
\]
\item[(b)] for every $m$,
\[
\mr{H}^m\left(\mr{R}_{s,t}(\mr{C})\right)/H^0_{\mathfrak m}(\mr{H}^m\left(\mr{R}_{s,t}(\mr{C})\right))
\]
is reflexive;
\item[(c)]
\[
\mathfrak m\,H^0_{\mathfrak m}\left(\mr{H}^m\left(\mr{R}_{s,t}(\mr{C})\right)\right)=0
\]
for every $m$;
\item[(d)] Conditions (a)--(c) above hold and the remaining Postnikov classes vanish.
\end{enumerate}

\item Conditions (a)--(c) above hold and the hyper-Tor spectral
sequence
\[
\operatorname{Tor}^R_p\left(\mr{H}^q\left(\mr{R}_{s,t}(\mr{C})\right),k\right)
\Longrightarrow
\mr{H}^{q-p}\left(\left(\mr{R}_{s,t}(\mr{C})\right)\otimes_R^Lk\right)
\]
degenerates at $E_2$.
\end{enumerate}
\end{theorem}
\2
\begin{remark}
In more geometric terms, with notation as above, the $\partial\bar\partial$-lemma is equivalent to requiring:
\begin{enumerate}
\item the cohomology sheaves $\widetilde{\mr{M}^m}$ be vector bundles over $\mr{U}$;
\item the sheaves $\widetilde{\mr{M}^m}/\widetilde{\mr{T}^m}$ satisfy the Hartogs' property at $0$.
\item The torsion subsheaf $\widetilde{\mr{T}^m}$ is scheme-theoretically supported at the reduced origin $\{0\}_{\text{red}}$.
\item The singularities at the origin are split, i.e. the Postnikov extensions of free and skyscraper sheaves vanish.
\end{enumerate}
\end{remark}
\subsection{\texorpdfstring{$\mr{E}_r$-}-degeneration and reduced torsion}
$\hspace*{3mm}$We conclude the section with a result relating degeneration at $\mr{E}_2$ and torsion analogously to the one-parameter torsion/flatness criteria for filtered complexes. 
\2
\begin{theorem}
Let $(\mr{C},\partial,\bar\partial)$ be a bounded double complex of finite dimensional vector spaces and let $$\mr{E}_2 = \mr{H}_\partial\left(\mr{H}_{\bar\partial}\right) \Rightarrow \mr{H}^{p+q}(\mr{Tot}(\mr{C}))$$ be the spectral sequence where we take $\bar\partial$-cohomology first. It collapses at $\mr{E}_r$ if and only if $$s^{r-1}H_s^0\left(\mr{H}^n(\mr{R}_{s,t}(\mr{C}))[t^{-1}]\right)=0$$ for all $n$, where $$H_s^0(\mr{M}) =\{x \in \mr{M} \; : \; s^Nx=0 \text{ for some } N\ge 1\}.$$
\end{theorem}
\2
\begin{remark}
Swapping the roles of $s$ and $t$, we obtain the same result for other the spectral sequence $$'\mr{E}_2 = \mr{H}_{\bar\partial}\left(\mr{H}_{\partial}\right) \Rightarrow \mr{H}^{p+q}(\mr{Tot}(\mr{C})).$$
\end{remark}
 \bibliography{references}
 \bibliographystyle{alpha}
\end{document}